\documentclass[11pt]{article}
\usepackage{amsmath}
\usepackage{amssymb}
\usepackage{amsbsy}
\usepackage{amsthm}
\usepackage{epsfig}
\usepackage{wrapfig}
\usepackage{eepic}
\usepackage{color}
\usepackage{graphicx}
\usepackage{amsfonts}%

\newtheorem{Proposition}{Proposition}[section]

\newtheorem{remark}{Remark}[section]
\newcommand{\bu}{{\bf u}}

\newcommand{\be}{{\bf e}}
\newcommand{\bD}{{\bf D}}

\newcommand{\bv}{{\bf v}}
\newcommand{\bw}{{\bf w}}
\newcommand{\bx}{{\bf x}}
\newcommand{\bn}{{\bf n}}

\newcommand{\bH}{{\bf H}}
\newcommand{\bX}{{\bf X}}

\newcommand{\bphi}{{\boldsymbol \phi}}

\def\div{\operatorname{div}}

\def\curl{\operatorname{curl}}
\newcommand{\blf} {\mathbf{f}}

\begin{document}

\title{On conservation laws of Navier-Stokes Galerkin discretizations}

\author{
Sergey Charnyi
\footnote{Department of Mathematical Sciences, Clemson University, Clemson, SC 29634 (scharny@clemson.edu).}
\and
Timo Heister
\footnote{Department of Mathematical Sciences, Clemson University, Clemson, SC 29634 (heister@clemson.edu), partially supported by NSF Grant DMS1522191
and by the Computational
Infrastructure for Geodynamics initiative (CIG), through NSF grant EAR-0949446.}
\and
Maxim A. Olshanskii
\footnote{Department of Mathematics, University of Houston, Houston TX 77004 (molshan@math.uh.edu), partially supported by Army Grant 65294-MA and NSF Grant DMS1522192.}
\and
Leo G. Rebholz\footnote{Department of Mathematical Sciences, Clemson University, Clemson, SC 29634 (rebholz@clemson.edu), partially supported by Army Grant 65294-MA and NSF Grant DMS1522191.}
}
\date{}

\maketitle

\begin{abstract}
We study conservation properties of  Galerkin methods for the incompressible Navier-Stokes equations, without the divergence constraint strongly enforced.  In typical discretizations such as the mixed finite element method, the conservation of mass is enforced only weakly, and this leads to discrete solutions which may not conserve energy, momentum, angular momentum, helicity, or vorticity, even though the physics of the Navier-Stokes equations dictate that they should.  We aim in this work to construct discrete formulations that conserve as many physical laws as possible without utilizing a strong enforcement of the divergence constraint, and doing so leads us to a new formulation that conserves each of energy, momentum, angular momentum, enstrophy in 2D, helicity and vorticity (for reference, the usual convective formulation does not conserve most of these quantities).  Several numerical experiments are performed, which verify the theory and test the new formulation.
\end{abstract}

\section{Introduction}

We consider formulations of the incompressible Navier-Stokes equations (NSE), which are given in a domain $\Omega\subset \mathbb{R}^d$, $d$=2 or 3, {\color{black}and for $t>0$} by
\begin{eqnarray}
\bu_t + (\bu\cdot\nabla) \bu + \nabla p - \nu\Delta \bu & = & \blf, \\
\div\bu & = & 0,\\
\bu(0) &= &\bu_0,
\end{eqnarray}
where $\bu$ and $p$ represent velocity and pressure, $\blf$ is an external forcing, $\bu_0$ is the initial velocity, and $\nu$ is the kinematic viscosity.  To solve this system, it must also be equipped with boundary conditions, {\color{black} and herein we consider no-slip, no penetration for velocity on the boundary: $\bu |_{\partial\Omega}={\bf 0}$.}  The NSE model the evolution of water, oil, and air flow (air under 220 m.p.h.), and therefore the ability to solve them is important in a wide array of engineering design problems.  Finding analytical solutions is known to be extremely difficult, and thus practitioners instead typically approximate solutions using numerical methods.

The purpose of this paper is to study conservation laws of solutions arising from discretizations of the NSE by Galerkin methods such as finite element methods, isogeometric analysis, or spectral element methods.  In typical discretizations, e.g., with Taylor-Hood finite elements, the conservation of mass is only weakly enforced, leading to discrete solutions $\bu_h$ which have $\div \bu_h \ne 0$.  While convergence of the $H^1$ error can often be proven, leading to the bound $\| \div \bu_h \| \le Ch^2$, the divergence error can still be significant on practical meshes (here $h$ is a characteristic step of an underlying mesh; in practice, there is a minimum $h$ that can produce solutions in reasonable time).  With the loss of mass conservation, it turns out that many other important conservation laws are also lost, including energy, momentum, angular momentum, and others, if steps are not taken in the development of the numerical discretization to make sure these quantities are conserved.  The fact that energy conservation is lost in Galerkin discretizations of the NSE is well-known, and a fix for this problem by using the skew-symmetric or rotation forms of the nonlinearity has been known for many years \cite{temam}.  A finite element  formulation for energy and helicity conservation was proposed in \cite{R07}, and in \cite{OR10b} it was discussed how an alternate (but equally valid) definition of helicity could be conserved by skew-symmetric formulations.  Similar phenomena happen with other types of discretization methods, and some clever discretizations have been developed which `bring back' conservation laws lost in standard discretizations, beginning decades ago by Arakawa, Fix, and others, for NSE and related equations \cite{AM03,A66,AL81,evans2013isogeometric,F75,LW04,PG16,ST89,SCN15}.  A common theme for all  `enhanced-physics' based schemes is that {\it the more physics is built into the discretization, the more accurate and stable the discrete solutions are, especially over longer time intervals}.

In the present work, we aim to develop numerical schemes/formulations that preserve even more conservation laws for the full NSE, beyond just energy.  By noticing that the key to discrete conservation properties is the formulation of the nonlinear term, we are able to find a  formulation of the NSE seemingly unconsidered in the literature, which conserves all of energy, momentum, and angular momentum; we call it the {\it Energy, Momentum, and Angular Momentum (EMA) conserving} formulation.  We propose this formulation in section 3, and formally show these conservation properties hold.  We also show that the usual convective, skew-symmetric and rotational formulations all fail to conserve momentum and angular momentum.  {\color{black} Additionally}, we show that the {\it EMA-conserving} formulation also conserves suitable definitions of vorticity, helicity, and 2D enstrophy.   Of course, if a Galerkin solution happens to be pointwise divergence-free, then all of the formulations are equivalent and each of them would conserve all of these quantities in an appropriate sense.  However, the use of such special element choices that provide pointwise divergence-free solutions (e.g. \cite{arnold:qin:scott:vogelius:2D,FN13,GN14,GN14b,Z05}) is not widespread, as they require constraints on the mesh and approximating polynomial degrees, and are not typically available in open source FE software for large scale computing \cite{dealII84}.

This paper is arranged as follows.  Section 2 presents notation and mathematical preliminaries that will allow for smoother analysis in later sections.  Section 3 presents the {\it EMA-conserving} formulation, and studies its conservation properties along with those of the convective, skew-symmetric, conservative, and rotational formulations.  Section 4 performs several numerical experiments, which test the conservation properties and accuracy of the various schemes.

\section{Notation and preliminaries}

Consider the domain $\Omega \subset \mathbb{R}^d$, {\color{black}$d$}=2 or 3,  and denote $(\cdot,\cdot)$ and $\| \cdot \|$ to be the $L^2(\Omega)$ inner product and norm on $\Omega$.

Consider $\bu,\ \bv,\ \bw\in \bH^1(\Omega)$,  and note that we do not enforce that any of these quantities are solenoidal except for the last two equations of this section.
Define the trilinear form $b:\bH^1(\Omega) \times \bH^1(\Omega) \times \bH^1(\Omega) \rightarrow \mathbb{R}$ by
\begin{equation}\label{conv_b}
b(\bu,\bv,\bw) = (\bu\cdot\nabla \bv,\bw).
\end{equation}
We recall the following properties of $b$.  The first two follow immediately from integration by parts, provided $\bu \in H_0^1(\Omega)$:
\begin{align}
b(\bu,\bv,\bw)& = -b(\bu,\bw,\bv)-( (\div \bu)\bv,\bw),  \label{vecid1} \\
b(\bu,\bw,\bw)& = -\frac12 \left( (\div \bu)\bw,\bw \right), \label{vecid2}\\
b(\bu,\bv,\bw)& = ( (\nabla \bv) \bu,\bw)=( (\nabla \bv)^T \bw,\bu).  \label{vecid3}
\end{align}
We denote the symmetric part of $\nabla \bu$ by
$
\nabla_s \bu :=  \bD(\bu) = \frac{ \nabla \bu + (\nabla \bu)^T }{2},
$
and the skew-symmetric part by
$
\nabla_n \bu := \frac{ \nabla \bu - (\nabla \bu)^T }{2}.
$
For any $\bu, \bv \in \bH^1(\Omega)$ one readily checks {\color{black} that}
\begin{equation}
(\nabla_n \bu)\bv =  \frac12 ( \curl{} \bu) \times \bv. \label{vecid6}
\end{equation}
Note that we define $\curl{} \bu$ in 2d in the usual way, as the 3d curl of $\bu$ extended by 0 in the
third component.

Straight-forward calculations provide the following  vector identities for functions $\bu, \bv \in \bH^1(\Omega)$:
\begin{align}
(\bu\cdot \nabla) \bu & =  (\curl \bu) \times  \bu + \nabla \frac12 | \bu |^2 =: (\curl \bu) \times  \bu + \nabla q, \label{vecid5} \\
(\nabla \bu)\bu & =  (\nabla_s \bu)\bu + (\nabla_n \bu)\bu = \bD(\bu)\bu +  \frac12 (\curl \bu) \times \bu, \label{vecid6b}
\end{align}
where $q:= \frac{ | \bu |^2}{2}$.
Also note that identity \eqref{vecid6b} implies that
\begin{equation}
( \bD(\bu) \bu,\bu) = (( \nabla \bu)\bu,\bu) = b(\bu,\bu,\bu).\label{vecid6c}
\end{equation}
From \eqref{vecid6}--\eqref{vecid6b} we obtain the following representation of the inertia term from the momentum equations:
\begin{equation}
(\bu\cdot \nabla) \bu = 2\bD(\bu)\bu - \nabla q. \label{vecid7}
\end{equation}
The identity \eqref{vecid7} is key to the new formulation we propose in the next section, which leads to improved discrete conservation properties.

\subsection{A parameterized vorticity equation}

Our study in section 3 of conservation laws {\color{black}for} vorticity and helicity involves different formulations of the vorticity equation.  We derive now a parameterized vorticity equation, which provides a family of {\color{black}formulations} which are equivalent when the velocity and vorticity are divergence-free.  {\color{black} Below, we denote vorticity by $\bw:=\curl \bu$.}


From \eqref{vecid6} with $\bv=\bw$, we find that $(\nabla_n \bu)\bw=0$, which implies
$(\nabla \bu)\bw = (\nabla \bu)^T \bw$.  The vorticity stretching term  $(\bw\cdot \nabla) \bu$ can thus be written as
\begin{eqnarray}
(\bw\cdot \nabla) \bu =  (\nabla \bu)\bw =  (\nabla \bu)^T\bw.
\end{eqnarray}
Similarly, the gradient of the helical density $\nabla (\bu \cdot \bw)$ can be written as
\[
\nabla (\bu \cdot \bw) = (\nabla \bu)^T\bw + (\nabla \bw)^T\bu = (\bw\cdot \nabla) \bu + (\nabla \bw)^T\bu.
\]
Hence, we can write the combination of the vorticity transport and stretching, as it appears in the equation's nonlinearity, for any real parameter $\beta_2$ as
\begin{eqnarray}
&& \hspace{-1in} (\bu\cdot\nabla) \bw - (\bw\cdot\nabla)\bu \nonumber \\
& = & (\bu\cdot\nabla) \bw - (\nabla \bu) \bw \nonumber \\
& = & (\bu\cdot\nabla) \bw -  \beta_2 (\nabla \bu)^T \bw - (1-\beta_2)(\nabla \bu) \bw \nonumber \\
& = & (\bu\cdot\nabla) \bw -  \beta_2 \left( \nabla(\bu \cdot \bw) - (\nabla \bw)^T \bu \right) - (1-\beta_2)(\nabla \bu) \bw. \label{vecid12}
\end{eqnarray}
Again using that $(\nabla \bu)\bw=(\nabla \bu)^T\bw$, we can also write
\[
(\nabla \bu) \bw = \beta_1 (\nabla \bu) \bw + (1-\beta_1)(\nabla \bu)^T \bw,
\]
for any real number $\beta_1$.  Denoting $\eta:=\bu\cdot \bw$, and using this together with \eqref{vecid12} and that $\div \bu=\div \bw=0$ (we have not made any divergence-free assumptions up to now), we obtain
\begin{multline}
\label{w_gen}
(\bu\cdot\nabla)\bw -(\bw\cdot\nabla)\bu =
(\bu\cdot\nabla)\bw + \beta_2(\nabla^T\bw)\bu  -(1-\beta_2)\left(\beta_1(\nabla \bu)+ (1-\beta_1)\nabla^T\bu\right)\bw\\  -\beta_2\nabla\eta + \beta_3(\div\bu)\bw+\beta_4(\div\bw)\bu,
\end{multline}
for real parameters $\beta_i$, $i=1,2,3,4$. The identities above lead us to the following vorticity equation, where the particular form of the nonlinearity depends on the choice of parameters:
\begin{align}
\bw_t + (\bu\cdot\nabla)\bw + \beta_2(\nabla^T\bw)\bu  -(1-\beta_2)\left(\beta_1(\nabla \bu)+ (1-\beta_1)\nabla^T\bu\right)\bw& \label{nse3} \\  -\nabla\eta + \beta_3(\div\bu)\bw+\beta_4(\div\bw)\bu - \nu\Delta \bw  &=  \curl{}  \blf, \nonumber \\
\div \bw & =  0.\label{nse4}
\end{align}
Note that the gradient term in \eqref{nse3} is not scaled with $\beta_2$. Hence, for $\beta_2\neq0$ the variable $\eta$ has the physical meaning of the scaled helical density, while for $\beta_2=0$ it is a {\color{black} Lagrange} multiplier corresponding to the divergence
free constraint and can be non-zero in the discrete setting.
We shall use discrete vorticity in the definition of  certain conserved quantities.  We will choose parameters  $\beta_i$ in such a way that the discrete vorticity solving the discrete counterpart of \eqref{nse3}--\eqref{nse4} delivers some desired conservation properties.

\section{Conservation properties under $\div \bu\neq {\bf 0}$ and the EMA formulation for Navier-Stokes}\label{s_4}

We now consider subspaces $\bX \subset \left[\bH^1_0(\Omega)\right]^d$, $Q \subset L^2(\Omega)$ of \textit{finite dimensions}.
To be more specific, we further assume that $\bX$ and $Q$ are finite element velocity and pressure spaces corresponding to an admissible triangulation of $\Omega$.
For simplicity we assume $\bX$ and $Q$ satisfy inf-sup compatibility conditions~\cite{GR86}; non inf-sup stable pairs require stabilization terms that will affect conservation properties, and should be studied separately and on a case-by-case basis.  We note our analysis can be easily extended to other Galerkin methods.

In most common discretizations of Navier-Stokes and related systems, the divergence-free constraint $\div \bu=0$ is only weakly enforced.  What holds instead of the pointwise constraint is that a numerical solution $\bu $ from $\bX $ satisfies
\[
(\div \bu ,q )=0 \quad \forall~ q \in Q ,
\]
where $Q $ is \textit{a finite dimensional} pressure space, for example piecewise linears which are globally continuous.  Even though convergence theory of mixed finite element methods exists that guarantees $\| \div \bu \|$ converges to 0 with optimal rate, in practical computations the divergence error can be large due to the associated constants being larger than the minimum practical meshwidth \cite{CELR11}. Enlarging the pressure space $Q $ to ensure $\div\bX \subset Q $ is usually not possible, since it would violate (apart of a few exceptional cases) the inf-sup compatibility condition
and make the method numerically unstable.

We now consider conservation properties of several common NSE formulations, along with a new one based on the identity \eqref{vecid7}.  To this end, we write the NSE momentum equation in the generic form:
\begin{eqnarray}
\bu_t + {\color{black} NL}(\bu) +\nabla p - \nu\Delta \bu  =  \blf , \label{ds1}
\end{eqnarray}
with the nonlinear terms defined for each formulation by
\begin{eqnarray*}
convective & : & {\color{black} NL}_{conv}(\bu) =  \bu \cdot\nabla \bu  \\
skew-symmetric & : & {\color{black} NL}_{skew}(\bu) =  \bu \cdot\nabla \bu +  \frac12 (\div \bu)\bu \\
rotational & : & {\color{black} NL}_{rot}(\bu) =  (\curl{} \bu)\times \bu \\
conservative & : & {\color{black} NL}_{cons}(\bu) =  \nabla \cdot (\bu{\color{black}\otimes}\bu) = \bu \cdot\nabla \bu +  (\div \bu)\bu.
\end{eqnarray*}
The convective, skew-symmetric, and rotational forms above are commonly  used in computational fluid dynamics and numerical analysis of fluid equations, see, e.g., \cite{G90, temam}, with the convective form being, probably, the most frequent choice in computation practice.

We propose now a new formulation, which we will show conserves energy, momentum and angular momentum, as well as helicity, vorticity, and 2D enstrophy, which we call the
\textit{energy momentum and angular momentum (EMA) conserving form}.  It is based on the following choice:
\begin{eqnarray*}
EMA \ conserving & : & {\color{black} NL}_{emac}(\bu) =  2\bD(\bu)\bu + (\div \bu)\bu.
\end{eqnarray*}
We remark that if we did assume the divergence constraint $\div \bu=0$ holds pointwise, then all above formulations are equivalent; for the EMA conserving scheme, this follows from \eqref{vecid7}.

The Galerkin method corresponding to various forms of inertia term reads: Find $\{\bu,p\}\in\bX\times Q$ satisfying
\begin{equation}
\left(\frac{\partial \bu }{\partial t} + {\color{black} NL}(\bu ),\bv \right) -(p ,\div\bv )+(q ,\div\bu ) +
\nu(\nabla\bu ,\nabla\bv )  =  (\blf,\bv )  \label{ds1FE}
\end{equation}
for all $\bv \in\bX $, $q \in Q $.

 For both the {\it rotational} and {\it EMA-conserving} formulations, the pressure $p$ is modified and includes a velocity contribution.  {\color{black} For the {\it rotational} formulation, the modified pressure is the Bernoulli pressure $p_{rot} = p_{kin} + \frac12 |\bu|^2$, where $p_{kin}$ is the kinematic pressure.  For the {\it EMA-conserving} formulation, the pressure is modified in a similar way, but with a negative sign: $p_{emac} = p_{kin} - \frac12 |\bu|^2$.}  To our knowledge, the {\it EMA-conserving}  formulation has yet to be considered in the literature, and our motivation for using it comes from Proposition \ref{conslaws1} below, which says that of these five formulations, only the {\it EMA-conserving} formulation exactly conserves energy,  momentum and angular momentum when the divergence constraint is not strongly enforced. Furthermore,
Proposition~\ref{conslaws2} in section~\ref{s_42} shows that this new formulation also exactly conserves suitably defined helicity, vorticity and 2D enstrophy.

\subsection{Energy,  momentum and angular momentum}
We now prove a result regarding conservation laws for \eqref{ds1FE}.  Our interest first is the conservation of energy, momentum and angular momentum:
\begin{eqnarray*}
&\text{Kinetic energy}\quad & E=\frac12(\bu,\bu):= \frac12 \int_\Omega |\bu|^2\mbox{d}\bx;\\
&\text{Linear momentum}\quad& M:=\int_\Omega \bu\, \mbox{d}\bx;\\
&\text{Angular momentum}\quad& M_\bx:=\int_\Omega \bu\times\bx\, \mbox{d}\bx.
\end{eqnarray*}

Most  useful boundary conditions alter the balance of  these quantities, as they should in the
presence of walls and interfaces. Moreover, the numerical treatment of boundaries, e.g. by enforcing conditions strongly or in a weak form, also affect this balance.
In this study, we  isolate the affect of the  treatment of the nonlinearity on the quantities of interest from the contribution of the boundary conditions. For this reason, we assume in section~\ref{s_4} that \textit{the finite element solution $\bu$ and $p$ is supported in some subset $\widehat{\Omega}\varsubsetneq \Omega$ of the computational domain $\Omega$}, i.e., there is a strip $S=\Omega \setminus \widehat{\Omega}$ along $\partial\Omega$ where $\bu $ is zero.  The same is assumed for the source term $\blf$.   We note this implies there is a strip of elements along the boundary where $\bu$, $p$, and $\blf$ vanish.
The prototypical scenario  is  the evolution of an isolated vortex in a self-induced flow.


\begin{Proposition} \label{conslaws1}
Assuming $(\div \bu ,q )=0$ for all $q \in Q $, but $\div \bu \ne 0$, the {\it skew-symmetric}, {\it rotational}, and {\it EMA-conserving} formulations conserve kinetic energy (for $\nu=0$, $\blf={\bf 0}$), and only the {\it EMA-conserving} and {\it conservative} formulations conserve momentum (for $\blf$ with zero linear momentum) and angular momentum (for $\blf$ with zero angular momentum).  Hence, the {\it EMA-conserving} is the only one of the formulations that conserves all three quantities.
\end{Proposition}
We divide the proof of this proposition into several subsections.

\subsubsection{Kinetic energy}
For energy conservation, testing \eqref{ds1FE} with $\bv =\bu $, $q =p $  gives
\[
\frac12 \frac{d}{dt} \| \bu  \|^2 + (NL(\bu ),\bu ) + \nu \| \nabla \bu  \|^2 = (\blf,\bu ).
\]
Thus, kinetic energy will be preserved for $\nu=0$, $\blf={\bf 0}$ if $(NL(\bu ),\bu )=0$.  For the {\it skew-symmetric} formulation, we use \eqref{vecid2} to get
\[
(NL_{skew}(\bu ),\bu ) = b(\bu , \bu ,\bu ) + \frac12((\div \bu )\bu ,\bu ) = 0,
\]
and for the {\it rot} formulation we use that the cross of two vectors is perpendicular to each of them,
\[
(NL_{rot}(\bu ),\bu ) = ((\curl{}  \bu )\times \bu ,\bu )  = 0.
\]
For the {\it EMA-conserving} formulation, we use \eqref{vecid6c} and then \eqref{vecid2} to obtain
\[
(NL_{emac}(\bu ),\bu ) = 2(\bD(\bu )\bu ,\bu ) +  ((\div \bu )\bu ,\bu ) =   2b(\bu ,\bu ,\bu ) +  ((\div \bu )\bu ,\bu )=0.
\]
For the {\it convective} formulation, the nonlinear term does not vanish in general:
\[
(NL_{conv}(\bu ),\bu ) =  b(\bu , \bu ,\bu ) = -\frac12  ((\div \bu )\bu ,\bu ),
\]
and thus kinetic energy will not be typically conserved by the {\it convective} formulation whenever $\div\bu \ne 0$.
Lastly, for the {\it conservative} formulation, we use the same identity as in the {\it convective} case, and find that
\[
(NL_{cons}(\bu ),\bu ) =  b(\bu , \bu ,\bu ) +  ((\div \bu )\bu ,\bu ) =  \frac12  ((\div \bu )\bu ,\bu ),
\]
and thus this formulation will not conserve kinetic energy in general.

\subsubsection{Momentum}
Next, we consider momentum conservation in the formulations. We cannot test \eqref{ds1FE} with $\bv=\be_i$ since this function is not in $\bX$.
Thanks to the assumption that $\bu\neq {\bf 0}$ only in some strictly interior subdomain $\widehat{\Omega}$,
we can define
the restriction
$\chi ({\bf g})\in \bX$ of an arbitrary function $\bf g$
by setting  $\chi (\bf g)=\bf g$ in $\widehat{\Omega}$ and
$\chi ({\bf g})$ arbitrary defined on $S=\Omega\setminus \widehat{\Omega}$ to satisfy zero boundary conditions.
We test \eqref{ds1FE} with $\bv =\chi (\be_i)\in \bX $  and $q =0$, which gives
\[
\frac{d}{dt}(\bu ,\be_i) + ({\color{black}NL}(\bu ),\be_i) = (\blf,\be_i),
\]
because the solution is zero on $S$.
Thus, momentum conservation is obtained if $(\blf,\be_i)=0$ and $({\color{black}NL}(\bu ),\be_i) =0$.  Thus we consider the latter for the different formulations.  In the {\it convective} formulation, we use \eqref{vecid1} and that $\be_i$ is constant to find that
\[
({\color{black}NL}_{conv}(\bu ),\be_i)= b(\bu,\bu ,\be_i) = - b(\bu,\be_i,\bu ) - ( (\div \bu )\bu ,\be_i) =  - ( (\div \bu )\bu ,\be_i),
\]
and for the {\it skew-symmetric} form we get
\[
({\color{black}NL}_{skew}(\bu ),\be_i) = b(\bu,\bu ,\be_i) + \frac12((\div \bu )\bu ,\be_i) = -\frac12 ((\div \bu )\bu ,\be_i).
\]
For {\it rotational} form, we use the vector identity $\bu \cdot\nabla \bu  = (\curl{} \bu )\times \bu  + \frac12 \nabla |\bu |^2$ to obtain
\begin{multline*}
({\color{black}NL}_{rot}(\bu ),\be_i) = ((\curl{}  \bu )\times \bu ,\be_i) = b(\bu,\bu ,\be_i) - \frac12 (\nabla |\bu |^2,\be_i)\\ = (\bu \cdot\nabla \bu ,\be_i) = - ( (\div \bu )\bu ,\be_i).
\end{multline*}
For the {\it EMA-conserving} formulation, however, the nonlinear term does vanish.  By expanding the rate of deformation tensor and using
$(\bu \cdot\nabla \bu ,\be_i) =  - ( ({\color{black} \div \bu} )\bu ,\be_i)$ and then \eqref{vecid3}, we find that
\begin{eqnarray*}
({\color{black}NL}_{emac}(\bu ),\be_i)  & = &  2(\bD(\bu )\bu ,\be_i) + ((\div \bu )\bu ,\be_i) \\
 & = &  b(\bu,\bu ,\be_i) + b(\be_i,\bu,\bu) + ((\div \bu )\bu ,\be_i) \\
 & = &  b(\be_i,\bu,\bu) \\
 & = & 0
\end{eqnarray*}
since $\be_i$ is divergence-free.

The {\it conservative} form also conserves momentum, as using the same identity as in the {\it convective} case, we obtain
\[
({\color{black}NL}_{cons}(\bu ),\be_i)=b(\bu, \bu ,\be_i) + ( (\div \bu )\bu ,\be_i)  =  - ( (\div \bu )\bu ,\be_i)  + ( (\div \bu )\bu ,\be_i) =0.
\]

\subsubsection{Angular momentum}
We consider next angular momentum conservation in the formulations; that is, whether or not they conserve $(M_\bx)_i:=(\bu ,\bphi_i),~~\bphi_i:=\bx\times\be_i,~~i=1,2,3$.  Note that $\div \bphi_i=0$ and $\Delta \bphi_i={\bf 0}$. Setting $\bv =\chi (\bphi_i)$, $q =0$ in \eqref{ds1FE} gives
\[
\left(\frac{\partial\bu }{\partial t},\bphi_i \right) + ({\color{black}NL}(\bu ),\bphi_i) + \nu(\nabla \bu ,\nabla\bphi_i)=(\blf,\bphi_i).
\]
Whether angular momentum is conserved comes down, once again, to whether it is preserved by the nonlinear term, i.e.{\color{black},} whether or not $({\color{black}NL}(\bu ),\bphi_i)=0$.  For the {\it EMA-conserving} formulation, since $\div \bphi_i=0$ we have that
\begin{eqnarray*}
({\color{black}NL}_{emac}(\bu ),\bphi_i)& = & 2(\bD(\bu )\bu ,\bphi_i) + ((\div \bu )\bu ,\bphi_i)  \\
& = & b(\bu,\bu ,\bphi_i) + b(\bphi_i,\bu,\bu) + ((\div \bu )\bu ,\bphi_i)  \\
& = &   b(\bu ,\bu ,\bphi_i) + ((\div \bu )\bu ,\bphi_i) \\
& =& -b(\bu,\bphi_i,\bu),
\end{eqnarray*}
with the last step coming from \eqref{vecid1}.  From here, expanding out the terms immediately reveals that $b(\bu,\bphi_i,\bu)=0$, and thus the {\it EMA-conserving} formulation does conserve angular momentum.

Similarly for the {\it conservative} formulation,
\begin{eqnarray*}
({\color{black}NL}_{cons}(\bu ),\bphi_i)& = & b(\bu,\bu ,\bphi_i) + ((\div \bu )\bu ,\bphi_i)  \\
& =& -b(\bu,\bphi_i,\bu) \\
& = & 0.
\end{eqnarray*}

For the {\it convective} formulation, similar identities reveal
\[
({\color{black}NL}_{conv}(\bu ),\bphi_i)=b(\bu,\bu ,\bphi_i)  = - ((\div \bu)\bu,\bphi_i ) \ne 0
\]
in general, and for the {\it skew-symmetric} formulation we use these same identities to obtain
\[
({\color{black}NL}_{skew}(\bu ),\bphi_i)= b(\bu,\bu ,\bphi_i) + \frac12 ((\div \bu)\bu,\bphi_i ) = - \frac12 ((\div \bu)\bu,\bphi_i ),
\]
which will not be zero in general either.  For the {\it rotational} formulation, we again use the vector identity
$\bu \cdot\nabla \bu  = (\curl{} \bu )\times \bu  + \frac12 \nabla |\bu |^2$, which provides since $\div \bphi_i=0$,
\[
({\color{black}NL}_{rot}(\bu ),\bphi_i)= ((\curl{} \bu ) \times  \bu ,\bphi_i) = (\bu\cdot \nabla \bu,\bphi_i) = - ((\div \bu)\bu,\bphi_i ),
\]
which is the same as for the {\it convective} formulation.

\subsection{Vorticity, helicity and 2D enstrophy}\label{s_42}
Denote by $\bu^*$ the exact Navier-Stokes solution  and $\bw^*=\nabla\times\bu^*$.
The following quantities are conserved for $\nu=0$ and suitable assumptions on the right-hand side $\blf$ and boundary conditions:
\begin{eqnarray*}
&\text{Helicity} \quad&H=(\bu^*,\bw^*):=\int_\Omega \bu^*\cdot\bw^* \mbox{d}\bx;\\
&\text{2D Enstrophy} \quad&H_{2D}=\frac12(\bw^*,\bw^*):= \frac12\int_\Omega \bw^*\cdot\bw^* \mbox{d}\bx\quad \text{(for a 2D flow)};\\
&\text{Total vorticity} \quad&W_i=(\bw^*,\be_i):=\int_\Omega w^*_i \mbox{d}\bx,\quad i=1,\dots,d.
\end{eqnarray*}
One verifies that none of the finite element methods discussed above conserve helicity, enstrophy or vorticity, see~\cite{OR10b}. It was further noted in~\cite{OR10b} that a more suitable definition of these quantities is based on the finite element solution to the vorticity equation rather than the curl of finite element velocity. This discrete vorticity still depends on the computed velocity $\bu $, but more implicitly, through the equation coefficients. Further, varying the coefficients $\beta_i$ in \eqref{nse3} we find the form of the finite element vorticity equation such that the recovered solution delivers conservation laws. We note that \textit{$\curl{} \bu $ does not necessarily ensure conservation of helicity, enstrophy or total vorticity even if the discrete solution $\bu $ is pointwise divergence free}. Defining discrete counterparts of these conserved quantities with the help of companion discrete vorticity equation is appropriate also in this case, see  the analysis of divergence-conforming B-splines for the unsteady Navier--Stokes equations in \cite{evans2013isogeometric}.

Consider the Navier-Stokes vorticity equation, which is found by taking the curl of the NSE:
\begin{equation}\label{vort0}
\bw^*_t + (\bu^* \cdot\nabla)\bw^* - (\bw^*\cdot\nabla)\bu^*  - \nu\Delta \bw^*  =  \curl{}  \blf.
\end{equation}

We consider now alternative discrete formulations that are equivalent to \eqref{vort0} when the velocity and vorticity are divergence-free, but which differ in discretizations.  A parametrized vorticity equation is given in \eqref{nse3}-\eqref{nse4}, which allows for such alternatives by the choice of $(\beta_1,\ \beta_2,\ \beta_3,\ \beta_4)$.  In the following discrete formulations, the velocity field $\bu\in\bX$ is the finite element solution to \eqref{ds1FE} and not the true NSE velocity $\bu^*$.  Note that {\color{black} homogeneous} Dirichlet boundary conditions are not appropriate for $\bw $ for general flow. However, since we assume $\bu$ vanishes in a neighborhood of $\partial \Omega$, then we \textit{assume} vorticity is also zero on and near the boundary.

For our first formulation of interest, we set  $\beta_1=1,\ \beta_2=0,\ \beta_3=0,\ \beta_4=0$.  This leads to the finite element formulation: find  $\bw \in\bX $ and Lagrange multiplier $\eta \in Q $   solving
\begin{multline}
\left(\frac{\partial\bw }{\partial t},\bv \right) + b(\bu ,\bw ,\bv )
-  b(\bw ,\bu ,\bv )
 + \nu(\nabla\bw ,\nabla\bv )+(\eta ,\div\bv )-(q ,\div\bw )  =  (\curl{} \blf,\bv ) \label{vort1}
\end{multline}
for $\bv \in \bX$ and $q \in Q$.  Alternatively, if we set $\beta_1=1,\ \beta_2=0,\ \beta_3=1,\ \beta_4=-1$, we arrive at the finite element formulation: find  $\bw \in\bX $ and Lagrange multiplier $\eta \in Q $   solving
\begin{multline}
\left(\frac{\partial\bw }{\partial t},\bv \right)
+ b(\bu ,\bw ,\bv ) -  b(\bw ,\bu ,\bv )
+ ((\div\bu )\bw ,\bv ) - ( (\div\bw )\bu ,\bv ) \\
+ \nu(\nabla\bw ,\nabla\bv )+(\eta ,\div\bv )-(q ,\div\bw )  =  (\curl{} \blf,\bv ) \label{vort2}
\end{multline}
for $\bv \in \bX$ and $q \in Q$.

For 2D flows, we consider the reduction to two dimensions after choosing $\beta_1=1, \beta_2=0,\ \beta_3=\frac12,\ \beta_4=0$, which provides the discrete formulation: find $w \in X$ satisfying for all $v\in X$,
\begin{equation}\label{vort3}
(w_t,v) +  ((\bu \cdot\nabla)w,v) + \nu(\nabla w,\nabla v)  + \frac12((\div\bu )w,v)
=  (\curl{} \blf,v).
\end{equation}


\begin{Proposition} \label{conslaws2}
Assume $\bu\in \bX$ solves \eqref{ds1FE} with the {\it EMA-conserving} form $N(\bu)=N_{emac}(\bu)$, and $\bw_1\in \bX$, $\bw_2\in \bX$, $(0,0,w)^T\in \bX$ are finite element vorticity solutions to \eqref{vort1}, \eqref{vort2}, \eqref{vort3}, respectively. The {\it EMA-conserving} formulation also conserves helicity (for $\blf=0$, $\nu=0$), 2D enstrophy (for $\curl{} \blf= {\bf 0}$, $\nu=0$), and  total vorticity in the sense of the following preserved quantities: $H=(\bu,\bw_1)$, $H_{2D}=\frac12\|w\|^2$, and $W_i=(\bw_2,\be_i)$.
\end{Proposition}

\begin{remark}\rm
For the other NSE formulations, conserved invariants involving vorticity also can be suitably defined.
\end{remark}

\subsubsection{Vorticity}\label{s_vort}
We start with the conservation of total vorticity and set $\bv =\chi (\be_i)$, $q =0$ in \eqref{vort2}.  Integration by parts immediately shows that $(\curl{} \blf,\be_i)=(\blf,\curl{} \be_i)=0$, and since $\be_i$ is constant, the $\eta$ term and the viscous term also drop, leaving
\begin{equation}
\left(\frac{\partial\bw_2 }{\partial t},\be_i \right)
+ b(\bu ,\bw_2 ,\be_i ) -  b(\bw_2 ,\bu ,\be_i )
+ ((\div\bu )\bw_2 ,\be_i ) - ( (\div\bw_2 )\bu ,\be_i ) =0.  \label{vort1n}
\end{equation}
 Since $\be_i$ is constant (and thus divergence-free), we have that $b(\bu ,\bw_2 ,\be_i) = - \left( (\div \bu )\bw_2 ,\be_i \right)$, and $b(\bw_2 ,\bu ,\be_i) = - \left( (\div \bw_2 )\bu ,\be_i \right)$, which reduces \eqref{vort1n} immediately to
 \[
\left(\frac{\partial\bw_2 }{\partial t},\be_i \right) = 0,
\]
and thus the total vorticity $W_i$ will be conserved for i=1,2,3.

\subsubsection{2D Enstrophy}\label{s_ens} For enstrophy conservation, take $v=w$ in \eqref{vort3} and set $\curl{} f={\bf 0}$ and $\nu=0$, which provides the equation
\begin{equation}
\frac12\frac{d}{dt} \| w \|^2 + b(\bu ,w,w) + \frac12((\div\bu )w,w) =0.
\end{equation}
Since $b(\bu ,w,w)=-\frac12 ((\div\bu )w,w)$, we have that $H_{2D}$ is conserved.

\subsubsection{Helicity}\label{s el}
For the conservation of $H$, set   $\bv =\bw_1 $, $q =0$ in \eqref{ds1FE} with $N=N_{emac}$, where $\bw_1$ solves \eqref{vort1}.   This vanishes the pressure term, and setting $\blf = {\bf 0}$ and $\nu=0$ yields
\begin{equation}
\left(\frac{\partial\bu }{\partial t},\bw_1 \right)  + 2(\bD(\bu )\bu ,\bw_1  )  + \left( (\div\bu )\bu ,\bw_1  \right)  =0. \label{hel1}
 \end{equation}
 Since $(\bD(\bu )\bu ,\bw_1  ) = \frac12 b(\bu ,\bu ,\bw_1 ) + \frac12 b(\bw_1 ,\bu ,\bu )$, we write \eqref{hel1} as
\begin{equation}
\left(\frac{\partial\bu }{\partial t},\bw_1 \right)  + b(\bu ,\bu ,\bw_1 ) +  b(\bw_1 ,\bu ,\bu )  + \left( (\div\bu )\bu ,\bw_1  \right)  =0. \label{hel2}
 \end{equation}
Next, take  $\bv =\bu$, $q =0$ in \eqref{vort1}, and with ${\bf f} = {\bf 0}$ and $\nu=0$ this provides
\begin{equation}
\left(\frac{\partial\bw_1 }{\partial t},\bv \right) + b(\bu ,\bw_1 ,\bu )
-  b(\bw_1 ,\bu ,\bu )   =  0. \label{hel3}
\end{equation}
Adding \eqref{hel2} and \eqref{hel3} gives the equation
\begin{equation}
\frac{d}{dt}(\bu,\bw_1) + b(\bu ,\bu ,\bw_1 )+ \left( (\div\bu )\bu ,\bw_1  \right)  + b(\bu ,\bw_1 ,\bu ) =0. \label{hel4}
\end{equation}
Using vector identity \eqref{vecid1}, we have that $b(\bu ,\bu ,\bw_1 )=-b(\bu ,\bw_1,\bu ) - ((\div \bu)\bu,\bw_1)$, which from
 \eqref{hel4} implies that $(\bu,\bw_1)$ is conserved.

{\color{black}
\subsection{Conservation properties of EMA-conserving formulation with Crank-Nicolson timestepping}
The analysis for the discrete conservation laws and balances performed above is for the semi-discrete case, i.e. without a temporal discretization. If a temporal discretization is applied, these balances can potentially be altered, as for example backward Euler and BDF2 are known to dissipate kinetic energy by their treatment of the time derivative terms.  We consider here the EMA-conserving Galerkin formulation together with Crank-Nicolson timestepping.  The Crank-Nicolson scheme is known to be energy conserving, and so seems a natural choice to study in this context.  We find that all the conservation properties that hold in the semidiscrete case also hold when Crank-Nicolson timestepping is used.

\begin{Proposition}
Solutions of the Crank-Nicolson  EMA-conserving  scheme \eqref{ds1FECN}-\eqref{ds2FECN} exactly conserve kinetic energy, momentum, angular momentum, helicity, vorticity, and 2D enstrophy (assuming no forcing or viscosity).
\end{Proposition}

Denoting half-steps of variables by $v^{n+\frac12}:=\frac12 \left( v^{n+1} +v^n \right)$, the scheme reads at each time step:  Given $\bu^n$ satisfying $(q ,\div\bu^{n} )  =  0$ for every $q\in Q$, find $\{\bu^{n+1},p^{n+\frac12} \}\in\bX\times Q$ satisfying
{\small
\begin{align}
 \left(\frac{\bu^{n+1} - \bu^n}{\Delta t},\bv\right)
+ \left( NL_{emac}(\bu^{n+\frac12}),\bv \right)
-(p^{n+\frac12} ,\div\bv )+
\nu(\nabla\bu^{n+\frac12} ,\nabla\bv )  & =   (\blf^{n+\frac12},\bv ),  \label{ds1FECN} \\
(q ,\div\bu^{n+1} ) & =  0,  \label{ds2FECN}
\end{align}}
for all $\bv \in\bX $, $q \in Q $.  Note that one solves directly for $p^{n+\frac12}$.
We will now analyze the conservation properties of \eqref{ds1FECN}-\eqref{ds2FECN}.  Just as above, we assume for the momentum and angular momentum that $\bu\ne {\bf 0}$ only in a strictly interior subdomain $\widehat{\Omega}\subset \Omega$.

As part of the the semi-discrete analysis above, we have already established that for any velocity field $\bu \in \bX$ it holds
\[
\left(NL_{emac}(\bu),\bu \right)   =
\left(NL_{emac}(\bu),\be_i \right) =
\left(NL_{emac}(\bu),\bphi_i \right)  =  0.
\]
Thus immediately we have that
\[
\left(NL_{emac}(\bu^{n+\frac12}),\bu^{n+\frac12} \right)    =
\left(NL_{emac}(\bu^{n+\frac12}),\be_i \right) =
\left(NL_{emac}(\bu^{n+\frac12}),\bphi_i \right)  =  0.
\]
Furthermore, since $\bu^{n+\frac12}$ is discretely divergence-free, i.e. $(q ,\div\bu^{n+1/2} )  =  0$ for every $q\in Q$,  and both $\be_i$ and $\bphi_i$ are divergence free, we also have that
\[
(p^{n+\frac12},\nabla \cdot \bu^{n+\frac12}) = (p^{n+\frac12},\nabla \cdot \be_i) = (p^{n+\frac12},\nabla \cdot \bphi_i) = 0.
\]
With these identities at hand, we proceed with the analysis, which shows that the Crank-Nicolson EMAC scheme conserves energy, momentum and angular momentum.

For conservation of energy, set $\bv = \bu^{n+\frac12}$, which from our identities above vanishes the nonlinear and pressure terms, leaving
\[
\frac{1}{2\Delta t} \left( \| \bu^{n+1} \|^2 - \| \bu^n \|^2 \right) + \nu \| \nabla \bu^{n+\frac12} \|^2 = (\blf^{n+\frac12},\bu^{n+\frac12}).
\]
Multiplying both sides by $\Delta t$ and summing over $M$ time steps provides
\[
\frac{1}{2} \| \bu^{M} \|^2 + \nu\Delta t\sum_{n=0}^{M-1} \| \nabla \bu^{n+\frac12} \|^2  = \frac{1}{2} \| \bu^{0} \|^2 + \Delta t \sum_{n=0}^{M-1} (\blf^{n+\frac12},\bu^{n+\frac12}),
\]
which is precisely the temporally discrete analog of the continuous-in-time energy balance (with integrals in time replaced by the composite midpoint approximation with rectangle width $\Delta t$).  Thus, if $\nu=0$ and $\blf={\bf 0}$, energy is exactly conserved.

For conservation of momentum, we take $\bv=\chi(\be_i)$, which vanishes the nonlinear and pressure terms, thanks to the identities above and the assumption that $\bu^{n+\frac12}$ vanishes on a strip along the boundary.  Since $(\nabla \bu^{n+\frac12},\nabla \chi(\be_i))= (\nabla \bu^{n+\frac12},\nabla \be_i)=0$, the viscous term also vanishes.  This leaves
\[
\frac{1}{\Delta t} \left( \bu^{n+1}  -  \bu^n , \be_i \right)  = (\blf^{n+\frac12},\chi(\be_i)),
\]
and so with the assumption that $\blf$ has zero momentum and vanishes outside of $\widehat{\Omega}$, exact momentum conservation is obtained:  for $i=1,2,3$,
\[
 \left( \bu^{n}   , \be_i \right)  =  \left(  \bu^0 , \be_i \right).
\]

Similarly for angular momentum, after choosing $\bv=\chi(\bphi_i)$, $\bphi_i=\bx\times\be_i$, $i=1,2,3$, immediately we get the nonlinear and pressure terms to vanish, and the assumptions on $\blf$ make the forcing term vanish as well.  This leaves
\[
\frac{1}{\Delta t} \left( \bu^{n+1}  -  \bu^n , \bphi_i \right)  + \nu (\nabla \bu^{n+\frac12},\nabla \bphi_i)=0.
\]
Hence if $\nu=0$, angular momentum is exactly conserved.

For helicity conservation, define $\bw^{n+1}$ to be the solution of: Given $\bw^n \in \bX$  satisfying $(q ,\div\bw )=0$ for every  $q\in Q$, and $\bu^{n+1},\bu^n$ solutions from the Crank-Nicolson EMAC scheme \eqref{ds1FECN}-\eqref{ds2FECN},
find $\bw^{n+1} \in\bX $ and Lagrange multiplier $\eta^{n+1} \in Q$  solving for $\bv \in \bX$ and $q \in Q$,
\begin{align}
\frac{1}{\Delta t} \left(\bw^{n+1} - \bw^n,\bv \right)
+ b(\bu^{n+\frac12} ,\bw^{n+\frac12} ,\bv ) -  b(\bw^{n+\frac12} ,\bu^{n+\frac12} ,\bv ) & \nonumber \\
+ \nu(\nabla\bw^{n+\frac12} ,\nabla\bv )+(\eta^{n+1} ,\div\bv )  & =  (\curl{} \blf^{n+\frac12},\bv ),  \label{dv11}  \\
(q ,\div\bw ) & =  0. \label{dv12}
\end{align}
Note that this vorticity need not be computed, and only its existence is required for the helicity conservation analysis.  Choose $\bv=\bu^{n+\frac12}$ in \eqref{dv11} and $\bv=\bw^{n+\frac12}$ in \eqref{ds1FECN}, and add these two equations.  Taking $\nu=0$, and $\blf={\bf 0}$, and noting that the pressure and $\eta$ terms drop,  we are left with
\begin{multline*}
\frac{1}{\Delta t} \left(\bw^{n+1} - \bw^n,\bu^{n+\frac12} \right)
+ \frac{1}{\Delta t} \left(\bu^{n+1} - \bu^n,\bw^{n+\frac12} \right)
+ b(\bu^{n+\frac12} ,\bw^{n+\frac12} ,\bu^{n+\frac12}) \\
-  b(\bw^{n+\frac12} ,\bu^{n+\frac12} ,\bu^{n+\frac12} )
 + (NL_{emac}(\bu^{n+\frac12})\bu^{n+\frac12},\bw^{n+\frac12}) = 0.
\end{multline*}
Now by identical arguments as the semi-discrete case (use $\bu^{n+\frac12}$ for $\bu$ and $\bw^{n+\frac12}$ for $\bw$), the nonlinear terms drop, leaving
\[
\frac{1}{\Delta t} \left(\bw^{n+1} - \bw^n,\bu^{n+\frac12} \right)
+ \frac{1}{\Delta t} \left(\bu^{n+1} - \bu^n,\bw^{n+\frac12} \right)
=0,
\]
which reduces to an exact conservation of helicity:
\[
\left(\bu^{n+1} ,\bw^{n+1} \right) = \left(\bu^{n} ,\bw^{n} \right)\quad\Rightarrow\quad
\left(\bu^{n} ,\bw^{n} \right) = \left(\bu^{0} ,\bw^{0} \right)\quad\forall\,n=1,2,\dots.
\]

For vorticity conservation, define $\bw^{n+1}$ to be the solution of: Given $\bw^0 \in \bX$  satisfying $(q ,\div\bw )=0$ for every  $q\in Q$, and $\bu^n$ ($n$=0,1,2,...,$n$+1) solutions from the Crank-Nicolson  EMA-conserving  scheme \eqref{ds1FECN}-\eqref{ds2FECN},
find $\bw^{n+1} \in\bX $ and Lagrange multiplier $\eta^{n+1} \in Q$  solving for $\bv \in \bX$ and $q \in Q$,
\begin{align*}
\frac{1}{\Delta t} \left(\bw^{n+1} - \bw^n,\bv \right)
+ b(\bu^{n+\frac12} ,\bw^{n+\frac12} ,\bv ) -  b(\bw^{n+\frac12} ,\bu^{n+\frac12} ,\bv ) & \\
+ ((\div\bu^{n+\frac12} )\bw^{n+\frac12} ,\bv ) - ( (\div\bw^{n+\frac12} )\bu^{n+\frac12} ,\bv )
+ \nu(\nabla\bw^{n+\frac12} ,\nabla\bv )& = (\curl{} \blf^{n+\frac12},\bv ).
\end{align*}
Similar to the helicity case, this vorticity also need not be computed, and only its existence is required for the vorticity conservation analysis. Taking  $\bv=\chi(\be_i)$, $\nu=0$, using the assumptions on $\blf$ and that $\be_i$ is constant in $\widehat{\Omega}$, we obtain
\begin{align*}
\frac{1}{\Delta t} \left(\bw^{n+1} - \bw^n,\be_i \right)
+ b(\bu^{n+\frac12} ,\bw^{n+\frac12} ,\be_i ) -  b(\bw^{n+\frac12} ,\bu^{n+\frac12} ,\be_i ) & \\
+ ((\div\bu^{n+\frac12} )\bw^{n+\frac12} ,\be_i ) - ( (\div\bw^{n+\frac12} )\bu^{n+\frac12} ,\be_i )
  & = 0.
  \end{align*}
Now by identical arguments as for vorticity conservation in the semi-discrete case, the nonlinear terms vanish, leaving exact momentum conservation: for $i$=1,2,3,
\[
(\bw^{n+1},\be_i)=(\bw^n,\be_i).
\]

Lastly, for 2D enstrophy conservation, let $w \in X$ satisfy for all $v\in X$,
\begin{eqnarray*}
\frac{1}{\Delta t}(w^{n+1}-w^n,v) +  ((\bu^{n+\frac12} \cdot\nabla)w^{n+\frac12},v) + \nu(\nabla w^{n+\frac12},\nabla v)  && \\
+ \frac12((\div\bu^{n+\frac12} )w^{n+\frac12},v)
&=&  (\curl{} \blf^{n+\frac12},v).
\end{eqnarray*}
The choice $v=w^{n+\frac12}$ immediately reveals conservation of 2D {\color{black}enstrophy}.

\subsection{Convergence of the EMA-conserving scheme}

Since the EMA-conserving scheme is new, it is important to check that it converges to the true solution, and with optimal rate.  Indeed, the analysis for the skew-symmetric case from \cite{Laytonbook} can be immediately extended to the EMA-conserving scheme, both for the semi-discrete case and the fully discrete case with Crank-Nicolson time stepping.  The key fact that allows the analysis to extend is the property of the nonlinearity that $(NL_{emac}(\bu),\bu)=0$.  Since this holds, the proofs from \cite{Laytonbook} of the skew-symmetric case can be completely mimicked.  For example, following the proof in \cite{Laytonbook} for the Crank-Nicolson, $(P_k,P_{k-1})$ velocity-pressure finite element discretization using the skew-symmetric form of the nonlinearity, one proves that under the usual assumptions,
\[
\max_{1\le n \le M} \| \bu^n - \bu_{nse}(t^n) \|
+ \left(\nu\Delta t \sum_{n=0}^{M-1} \| \nabla (\bu^{n+\frac12} - \bu_{nse}(t^{n+\frac12})) \|^2 \right)^{1/2}
\le
C(h^k + \Delta t^2),
\]
where $\bu_{nse}$ is a true solution to the Navier--Stokes equations.  Proofs for other time stepping methods for the Navier--Stokes equations with skew-symmetric formulation can also be immediately adapted for the EMA-conserving formulation.

}

\subsection{Discussion}
We have now established that the {\it EMA-conserving} formulation does indeed conserve energy, momentum, angular momentum, and suitable definitions of enstrophy (in 2D), helicity, and vorticity.  One may question if the {\it EMA-conserving} formulation is the only one or the `simplest' one which conserves all quantities listed above. We do not have an ultimate answer to these questions. Nevertheless, attempting to address it let us comment on the way we deduce this formulation: Similar to the vorticity equation \eqref{nse3}-\eqref{nse4}, we can write the momentum equation with linear combinations of different forms of the inertia terms from \eqref{conv_b}, \eqref{vecid5}, \eqref{vecid6b} and additional divergence terms. The {\it EMA-conserving} formulation is then found to be the \textit{unique} combination that conserves discrete kinetic energy, momentum, and angular momentum.  As already discussed, the  conservation of the discrete helicity, 2D enstrophy and vorticity are further understood with the help of the {\color{black}companion} finite element vorticity equations. {\color{black}We stress that the vorticity equation is \textit{not} a part of the finite element method here. It is introduced only to suitably define discrete conserved quantities. Nevertheless, the finite element vorticity equation can be used for postprocessing the finite element velocity in order to recover physically `correct' vorticity, if desired.}

\section{Numerical Experiments}

We now provide results of several numerical experiments that test and compare the different NSE formulations.  The specific formulations we test are (for the case of homogeneous Dirichlet boundary conditions): Find $(\bu_h,p_h) \in (X_h,Q_h)$ such that for every $(\bv_h,q_h)\in (X_h,Q_h)$, \\
Convective formulation (CONV)
\begin{eqnarray*}
((\bu_h)_t ,\bv_h) + (\bu_h\cdot\nabla \bu_h,\bv_h) - (p_h,\div  \bv_h) + \nu(\nabla \bu_h,\nabla \bv_h)  & = & (\blf, \bv_h), \\
(\div  \bu_h,q_h) & = & 0.
\end{eqnarray*}
Skew-symmetric formulation (SKEW)
\begin{eqnarray*}
((\bu_h)_t ,\bv_h) + (\bu_h\cdot\nabla \bu_h,\bv_h) + \frac12 ( (\div  \bu_h)\bu_h,\bv_h) - (p_h,\div  \bv_h) + \nu(\nabla \bu_h,\nabla \bv_h)  & = & (\blf, \bv_h), \\
(\div  \bu_h,q_h) & = & 0.
\end{eqnarray*}
Conservative formulation (CONS)
\begin{eqnarray*}
((\bu_h)_t ,\bv_h) + (\bu_h\cdot\nabla \bu_h,\bv_h) +  ( (\div  \bu_h)\bu_h,\bv_h) - (p_h,\div  \bv_h) + \nu(\nabla \bu_h,\nabla \bv_h)  & = & (\blf, \bv_h), \\
(\div  \bu_h,q_h) & = & 0.
\end{eqnarray*}
Rotational formulation (ROT)
\begin{eqnarray*}
((\bu_h)_t ,\bv_h) +((\curl{}  \bu_h)\times \bu_h,\bv) - (p_h,\div  \bv_h) + \nu(\nabla \bu_h,\nabla \bv_h)  & = & (\blf, \bv_h), \\
(\div  \bu_h,q_h) & = & 0.
\end{eqnarray*}
Energy, momentum, and angular momentum conserving formulation (EMAC)
\begin{eqnarray*}
((\bu_h)_t ,\bv_h) + 2(D(\bu_h)\bu_h,\bv_h) + ((\div  \bu_h)\bu_h,\bv_h) - (p_h,\div  \bv_h) + \nu(\nabla \bu_h,\nabla \bv_h)  & = & (\blf, \bv_h), \\
(\div  \bu_h,q_h) & = & 0.
\end{eqnarray*}

For the temporal discretizations, our tests employ several temporal discretizations, including Crank-Nicolson method for the Gresho problem described below (since here we test for integral invariants), BDF2, and BDF3.  The choice of Taylor-Hood velocity-pressure elements is used throughout, which is $(P_2,P_1)$ on triangular meshes, and $(Q_2,Q_1)$ on quadrilateral meshes. In the latter case, the deal.II finite element library (\cite{dealII84}) {\color{black} was} used.  No stabilization was used in any of the 2D simulations, however for the $(Q_2,Q_1)$ computations, grad-div stabilization~\cite{reusken} with a small parameter ($\gamma=0.1$) was used since it is an integral part of the preconditioner used for the linear solves.  We recognize that different element choices and different stabilizations can improve these schemes to varying degrees; future studies certainly could include various stabilization and element choices.

{\color{black}
The nonlinear problem  in each time step is solved using Newton's method
with a tolerance of $10^{-8}$, usually requiring only 1 or 2 iterations.
The corresponding linear systems  are solved in parallel using the grad-div based block preconditioner based on \cite{HR13}:
We solve the linear system $Mx=b$ with
\[ M = \begin{pmatrix} A & B^T \\ B & 0 \end{pmatrix}
\]
using FGMRES~\cite{saad1993flexible} with right preconditioning.  The  preconditioner {\color{black} $P^{-1}$} has the form
\[
 P^{-1} = \begin{pmatrix} \hat{A} & B^T \\ 0 & \hat{S} \end{pmatrix}^{{\color{black}-1}}.
\]
{\color{black}Here $\hat{A}^{-1}$ is defined through an inner iterative process with the matrix $A$}.  For this purpose we used  the GMRES method  preconditioned with an algebraic multigrid V-cycle and a relative tolerance of 1e-3.
The matrix $\hat{S}^{-1}$, {\color{black} approximating the inverse of the Schur complement $S$}, is defined by the sum
\[
 \hat{S}^{-1} = S_1^{-1} + S_2^{-1},
 \]
 where
 $ S_1 = 1/(\nu+\gamma) M_p $
 and
 $ S_2 = 1.0/c L_p $. Here, $M_p$ and $L_p$ are pressure mass and Laplace matrices and $1/c$ is the timestep size.
The action for $S_1^{-1}$ and $S_2^{-1}$ are approximated by separate
inner  GMRES solves preconditioned with block ILU(0).}

For the channel flow problems with an outflow, we weakly enforce the zero-traction boundary condition $  (-\nu\nabla \bu + pI)\cdot \bn |_{\Gamma_{out}} =0$.  For the CONV and CONS formulations, this becomes a `do-nothing' condition.  For the rest of the formulations, it requires a nonlinear boundary integral at the outflow.
\medskip

To illustrate the conservation properties of the various formulations, we choose several test problems:
For the first one, the quantities of interest are exactly conserved, while other test cases represent more realistic
scenarios of viscous fluid flows passing streamlined or bluff bodies. In the latter case, viscous and boundary effects perturb all conservation laws. We include these test cases in the attempt to give the  first assessment of other properties of the EMAC form such as numerical stability and accuracy.

\subsection{Gresho Problem}

\begin{figure}[!h]
\begin{center}
\includegraphics[width=.31\textwidth,height=0.31\textwidth, viewport=70 45 550 400, clip]{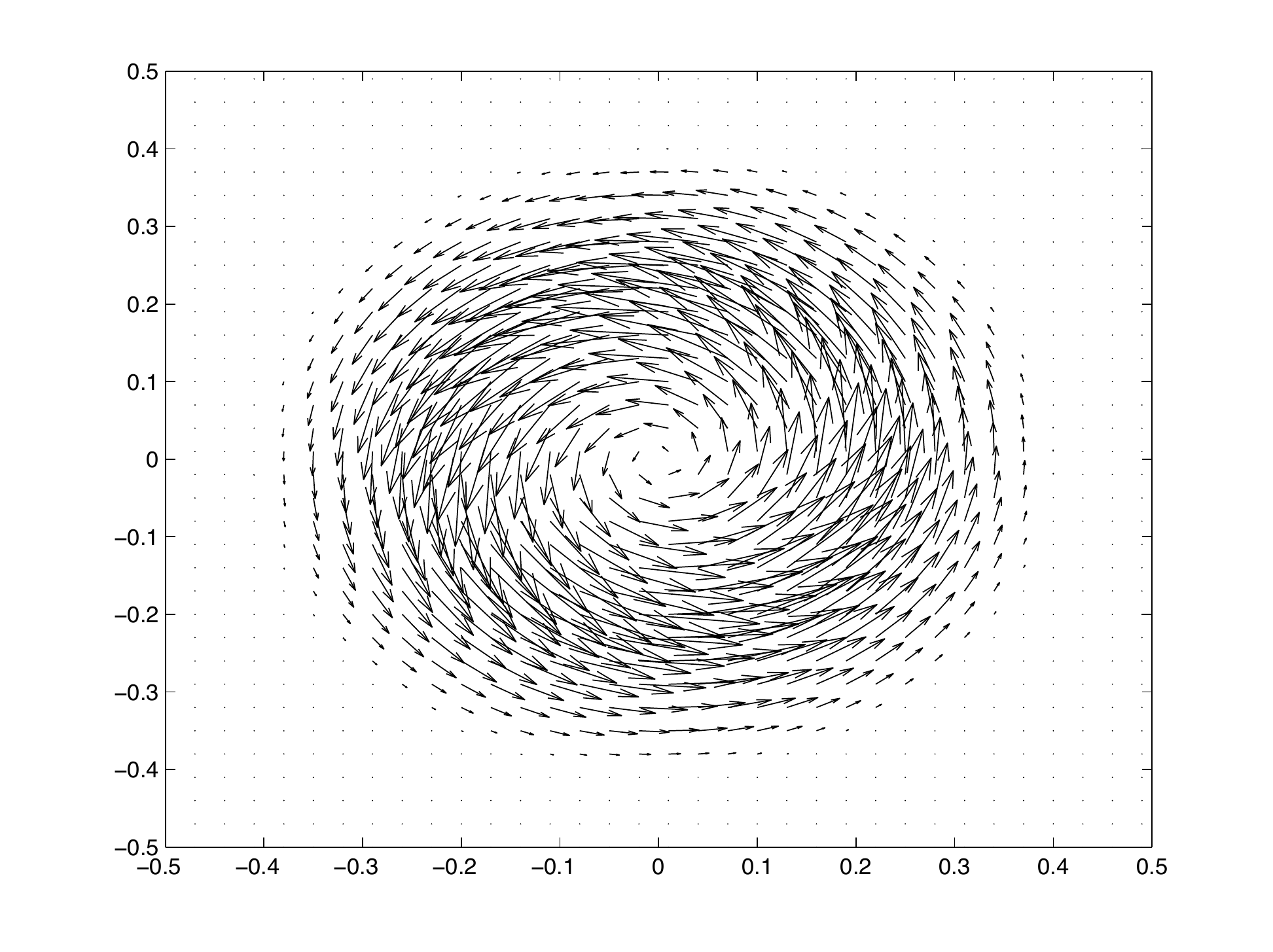}
\includegraphics[width=.35\textwidth,height=0.31\textwidth, viewport=70 45 550 440, clip]{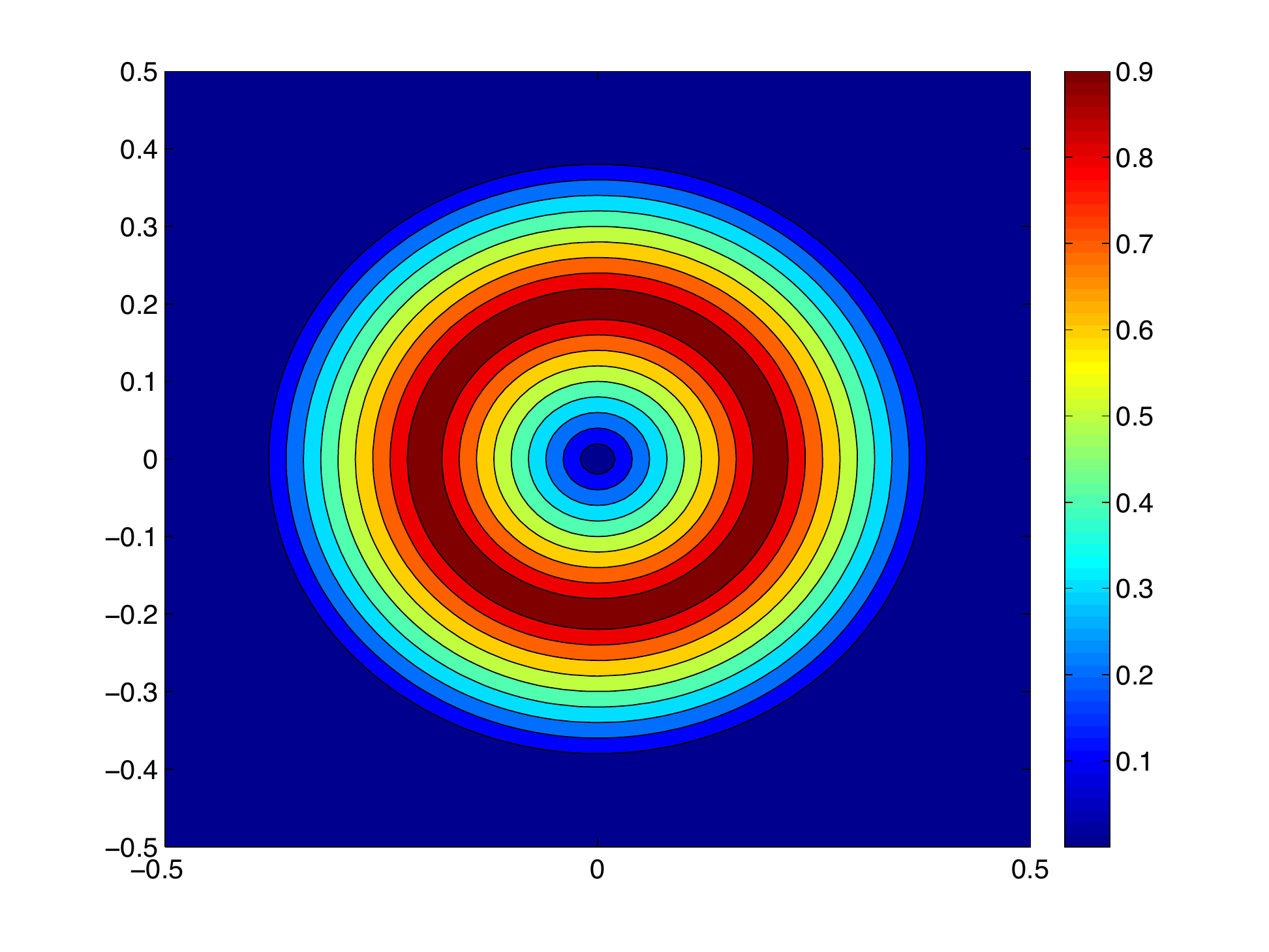}
\end{center}
\caption{\label{vortexmesh}
Shown above is the true velocity solution for the Gresho problem as a vector plot (left) and speed contour plot (right). }
\end{figure}

\begin{figure}[!h]
\begin{center}
\includegraphics[width=.49\textwidth,height=0.25\textwidth, viewport=60 0 860 320, clip]{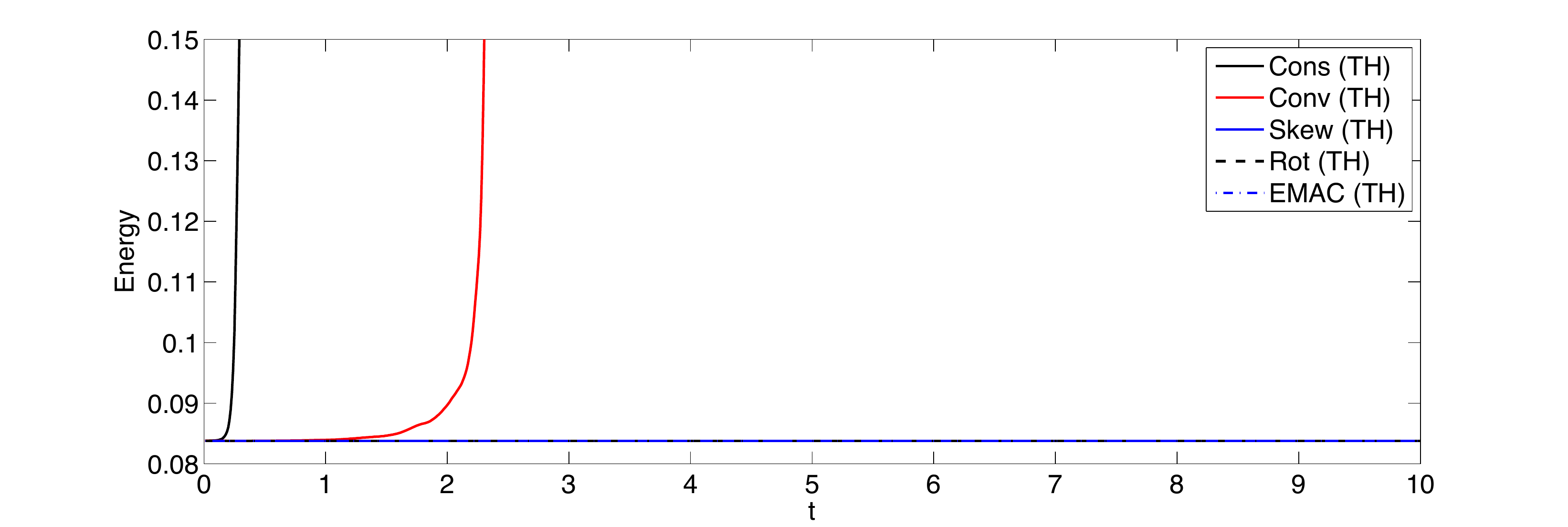}
\includegraphics[width=.49\textwidth,height=0.25\textwidth, viewport=60 0 860 320, clip]{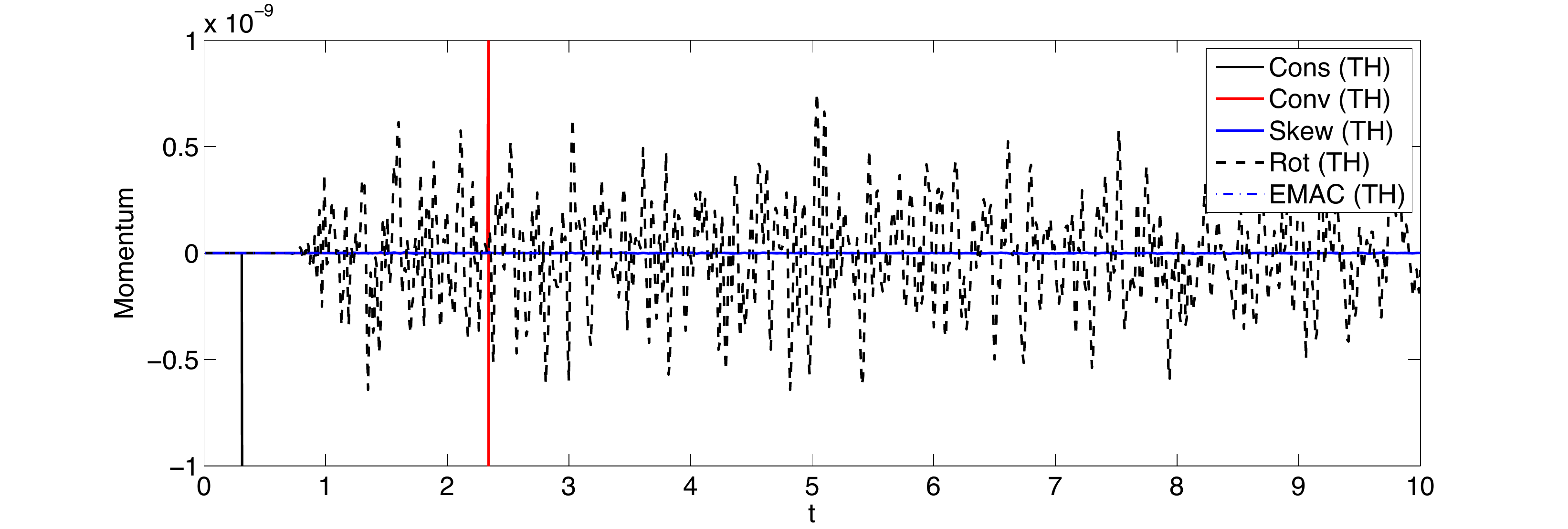}
\includegraphics[width=.49\textwidth,height=0.25\textwidth, viewport=60 0 860 320, clip]{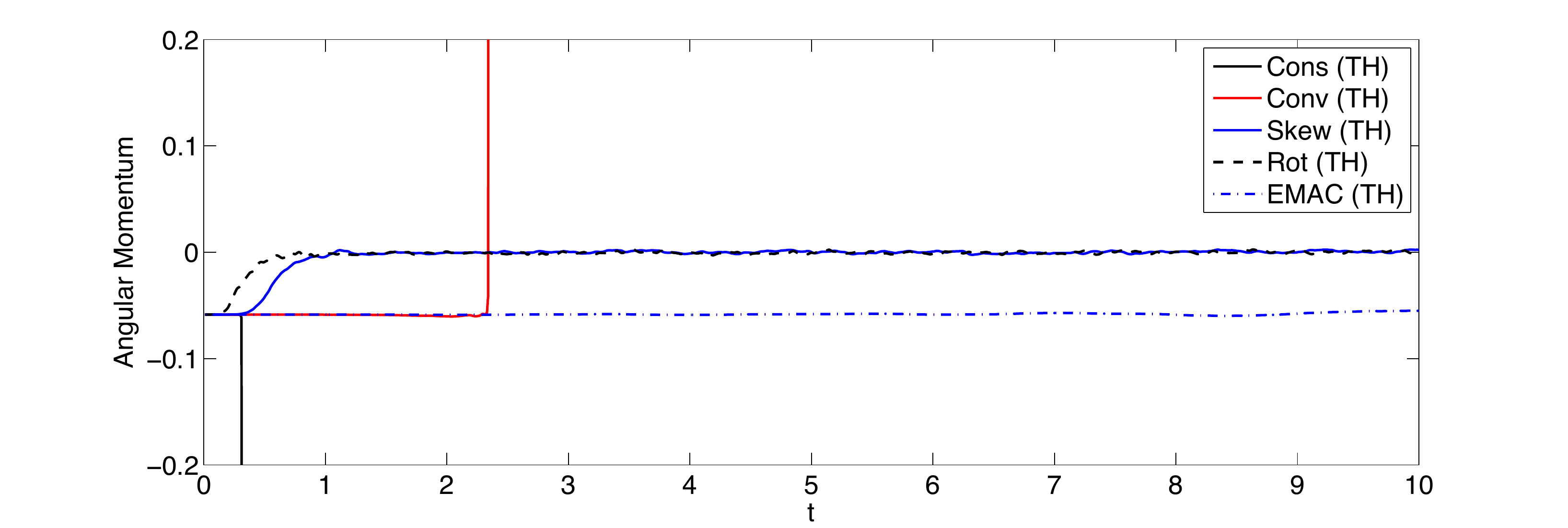}
\includegraphics[width=.49\textwidth,height=0.25\textwidth, viewport=60 0 860 320, clip]{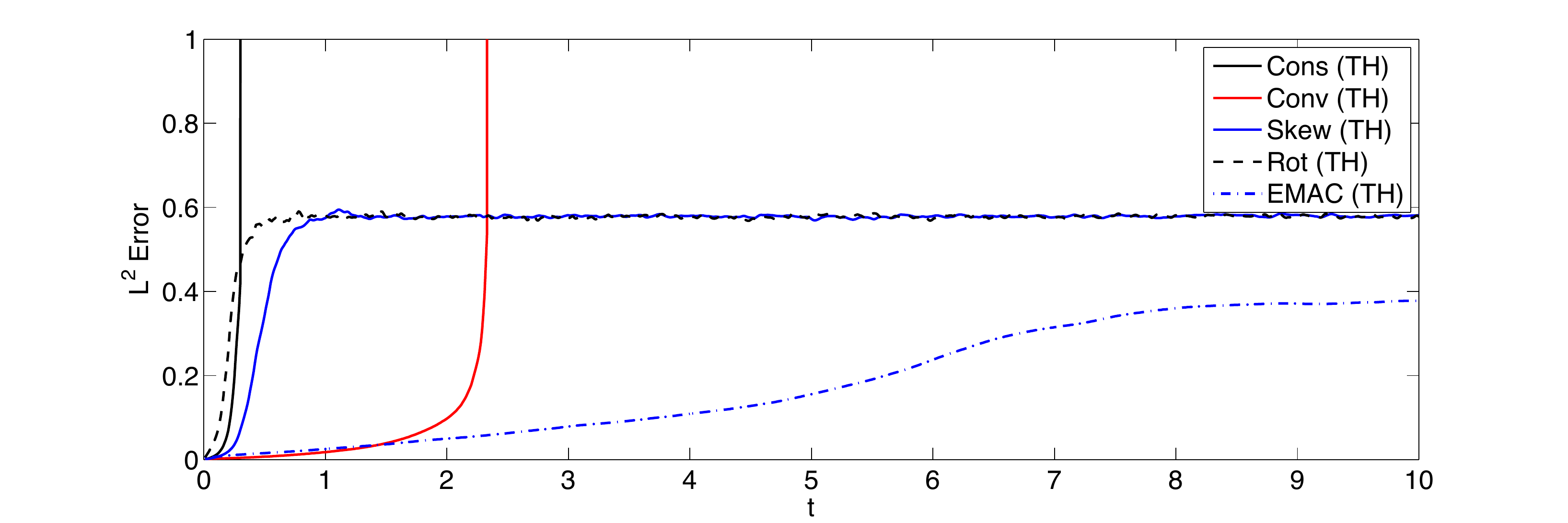}
\end{center}
\caption{\label{Gresho2}
Shown above are plots of time versus energy, momentum, angular momentum, and $L^2(\Omega)$ velocity error, for the various formulations in the Gresho problem, using Taylor-Hood elements.}
\end{figure}


We consider first the Gresho problem, which is often referred to as the `standing vortex problem' \cite{TMRS92,LW03,G90}.  The problem is defined by starting with an initial condition $\bu_0$ that is an exact solution of the steady Euler equations.  On $\Omega=(-.5,.5)^2$, with $r=\sqrt{x^2+y^2}$, the velocity and pressure solutions are defined by
\[
\begin{aligned}
r \le 0.2~: \left\{\begin{split} \bu &= \left( \begin{array}{c} -5y \\ 5x \end{array} \right) \\
 p&=12.5r^2 + C_1 \end{split}\right.,\quad
r>0.4 ~:  \left\{ \begin{split} \bu& = \left( \begin{array}{c} 0 \\ 0 \end{array} \right)\\ p&=0\end{split}\right.,  \\
0.2 \le r \le 0.4~: \left\{\begin{split}\bu &= \left( \begin{array}{c}  \frac{2y}{r}+5y \\ \frac{2x}{r} - 5x \end{array} \right)   \\
 p&=12.5r^2 - 20r + 4\log(r) + C_2 \end{split}\right. ,
\end{aligned}
\]
where
\[
C_2 = (-12.5)(0.4)^2 + 20(0.4)^2 - 4\log(0.4),\ C_1 = C_2 - 20(0.2) + 4\log(0.2).
\]
The vorticity ($w={u_2}_x - {u_1}_y$) can be calculated to be $w=10$ when $r\le 0.2$, $w=2/r-10$ on $0.2\le r\le 0.4$, and $w=0$ when $r>0.4$.  This is an interesting problem because it is an exact solution of the steady Euler equations, i.e.
\[
\bu\cdot\nabla \bu + \nabla p = 0.
\]
Since we choose the initial condition to be this steady Euler solution, an accurate scheme should preserve the solution in time.  Moreover, it is also a good test for a numerical scheme's ability to conserve certain quantities such as energy, momentum and angular momentum, since no viscosity or forcing is present, and the boundaries do not play a role (unless significant error causes nonzero velocity to creep out to the boundary).  A plot of the true velocity solution is shown in Figure \ref{vortexmesh}.

We compute solutions to the Gresho problem using the different formulations, together with Crank-Nicolson time stepping (using Newton's method to solve the nonlinear problem at each time step), with ${\bf f}={\bf 0}$, $\nu=0$, and no-penetration boundary conditions up to T=10.  We computed using $(P_2,P_1)$ Taylor-Hood elements on a 48x48 uniform mesh  and a time step of $\Delta t=0.01$.

{\color{black}
Plots of energy, momentum, angular momentum, and $L^2$ velocity error versus time are shown in Figure \ref{Gresho2}.  The EMAC scheme gives the best results: it conserves energy and momentum, is the only scheme to conserve angular momentum, and has significantly better $L^2(\Omega)$ error than all the other methods.   The CONS scheme gives by far the worst results.  The energy of the CONS solution is blowing up, which causes the nonlinear solver to fail before t=0.20.   CONV also blows up, although not until around t=2.4; until it blows up, it gives errors similar to EMAC.  ROT and SKEW are energy conserving and stable, but have poor accuracy compared to EMAC.  We note that all the results for conserved quantities are consistent with the theory of the previous section, and in particular the EMAC scheme is the only one to conserve each of energy, momentum and angular momentum.
}

\subsection{Channel flow around a cylinder}

\begin{figure}[h!]
\begin{center}
\includegraphics[width=0.7\textwidth,height=0.24\textwidth, trim=0 0 0 0, clip]{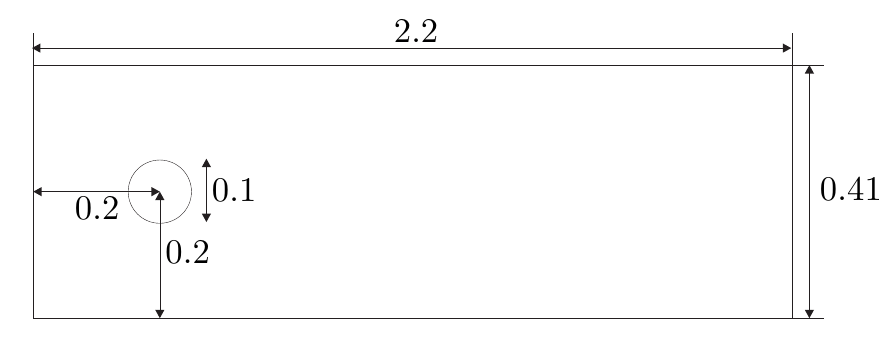}\\
\includegraphics[width=0.65\textwidth,height=0.17\textwidth, trim=180 330 70 30, clip]{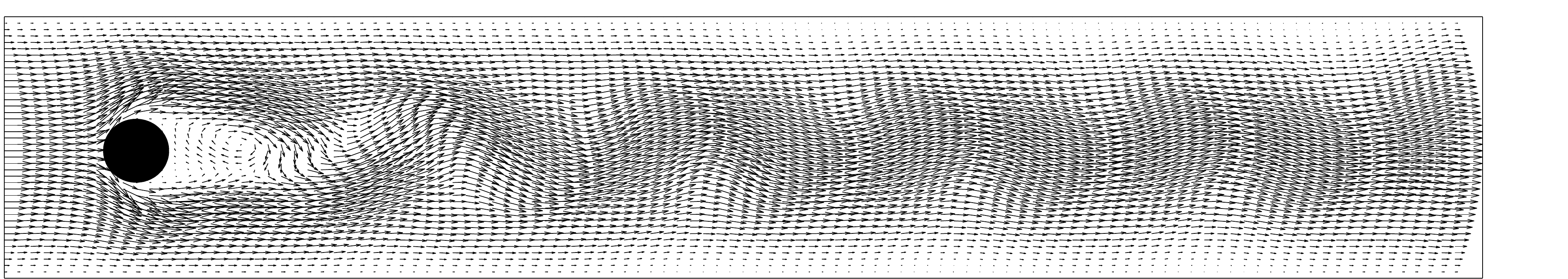}
\end{center}
\caption{\label{cyldomain} Shown above is the channel flow around a cylinder domain (top), and a resolved velocity field at t=6.}
\end{figure}

Our next experiment tests the algorithms above on the flow around
a cylinder benchmark problem, taken from \cite{J04,ST96}.  The domain for the problem is a $2.2\times 0.41$ rectangular channel with a circle (cylinder) of radius $0.05$ centered at $(0.2,0.2)$, see Figure \ref{cyldomain}.

No slip boundary conditions are enforced on the walls and cylinder, and
the time dependent inflow profile is taken to be
\begin{align*}
& u_{1}(0,y,t)=u_{1}(2.2,y,t)=\frac{6}{0.41^{2}}\sin(\pi\,t/8)y(0.41-y) \, ,
\\
& u_{2}(0,y,t)=u_{2}(2.2,y,t)=0,
\end{align*}
and a zero-traction outflow condition is weakly enforced.  The viscosity is set as $\nu=10^{-3}$ and there is no external force, ${\bf f}={\bf 0}$.

This problem is well studied, and it is known that as the flow rate increases, two vortices start to develop
by T=4 behind the cylinder. They then separate into the flow, and soon after a vortex street forms which can be
visible through t=8. Reference values for lift and
drag coefficients, and for pressure drop across the cylinder at t=8 are given in \cite{J04} as
\begin{equation*}
c_{d,max}^{ref}=2.95092 , \quad c_{l,max}^{ref}=0.47795, \quad \Delta p^{ref}=-0.11160.
\end{equation*}

\begin{table}[h!]
 \centering
 \begin{tabular}{|c|c|c|c|c|c|c|c|c|}  \hline
Method & $dim(X_h)$ & $\Delta t$   & $c_d^{max}$ & $|error|$ & $c_l^{max}$ & $|error|$ & $\Delta p$(8) & $|error|$  \\ \hline
ROT   & 34,762 & 0.005 & 2.94442 & 6.48E-3 & 0.412069 & 6.59E-2 & -0.11168 & 8.20E-5 \\ \hline
CONV  & 34,762 & 0.005 & 2.94672 & 4.18E-3 & 0.470062 & {\bf 7.94E-3} & -0.11176 & 1.62E-4 \\ \hline
SKEW  & 34,762 & 0.005 & 2.94678 & 4.12E-3 & 0.467538 & 1.05E-2 & -0.11177 & 1.70E-4 \\ \hline
CONS  & 34,762 & 0.005 & 2.94667 & 4.25E-3 & 0.450239 & 2.77E-2 & -0.11179 &  1.90E-4\\ \hline
EMAC   & 34,762 & 0.005 & 2.94819 & {\bf 2.71E-3} & 0.525675 & 4.77E-2 & -0.11166 & {\bf 5.68E-5} \\ \hline \hline
ROT  & 61,694 & 0.005 & 2.94638 & 4.52E-3 & 0.484486 & 6.49E-3 & -0.11139 & 2.10E-4 \\ \hline
CONV  & 61,694 & 0.005 & 2.94893 & 1.97E-3 & 0.478282 & {\bf 2.82E-4} & -0.11159 & {\bf 1.13E-5} \\ \hline
SKEW & 61,694 & 0.005 & 2.94892 & 1.98E-3 & 0.477249 & 7.51E-4 & -0.11158 & 2.15E-5 \\ \hline
CONS & 61,694 & 0.005 & 2.94891 & 1.99E-3 & 0.477013 & 9.37E-4 & -0.11149 & 1.10E-4 \\ \hline
EMAC   & 61,694 & 0.005 & 2.94961 & {\bf 1.29E-3} & 0.490655 & 1.27E-2 & -0.11119 & 4.06E-4 \\ \hline \hline
ROT   & 95,738 & 0.005 & 2.94919 & 1.71E-3 & 0.480021 & 2.02E-3 & -0.11186 &  2.64E-4 \\ \hline
CONV  & 95,738 & 0.005 & 2.94966 & 1.24E-3 & 0.478567 & 5.67E-4 & -0.11155 & {\bf 5.00E-5} \\ \hline
SKEW  & 95,738 & 0.005 & 2.94966 & 1.24E-3 & 0.478106 & {\bf 1.06E-4} & -0.11154 & 6.04E-5 \\ \hline
CONS  & 95,738 & 0.005 & 2.94966 & 1.24E-3 & 0.477831 & 1.19E-4 & -0.11155 & {\bf 5.00E-5} \\ \hline
EMAC   & 95,738 & 0.005 & 2.94986 & {\bf 1.04E-3} & 0.484425 & 6.43E-3 & -0.11141 & 1.93E-4 \\ \hline
\end{tabular}%
\caption{ \label{ldtable} Max lift and drag coefficients, and pressure drop across the cylinder at t=8, for the various formulations, using $(P_2,P_1)$ elements.}
\end{table}

We computed solutions using several meshes with Taylor-Hood elements, BDF3 time stepping, and time step $\Delta t=0.005$ (we also used $\Delta t=0.01$ and obtained very similar results).  Results for maximum lift and drag, as well as for the t=8 pressure drop are shown in Table \ref{ldtable}.  For each mesh, the best errors are made bold for each statistic.  We observe that in each case, the EMAC formulation provides the best prediction of the maximum drag coefficient, CONV and SKEW forms provide the best maximum lift coefficient prediction, and the EMAC, CONV and CONS provide the best predictions of pressure drop error.  Overall, the methods give rather similar predictions, and it is fair to say the methods are comparable for this test problem with these discretizations.

\subsection{Channel flow past a flat plate at Re=100}\label{s_exp}

\begin{figure}[ht]
\begin{center}
\includegraphics[width=0.75\textwidth]{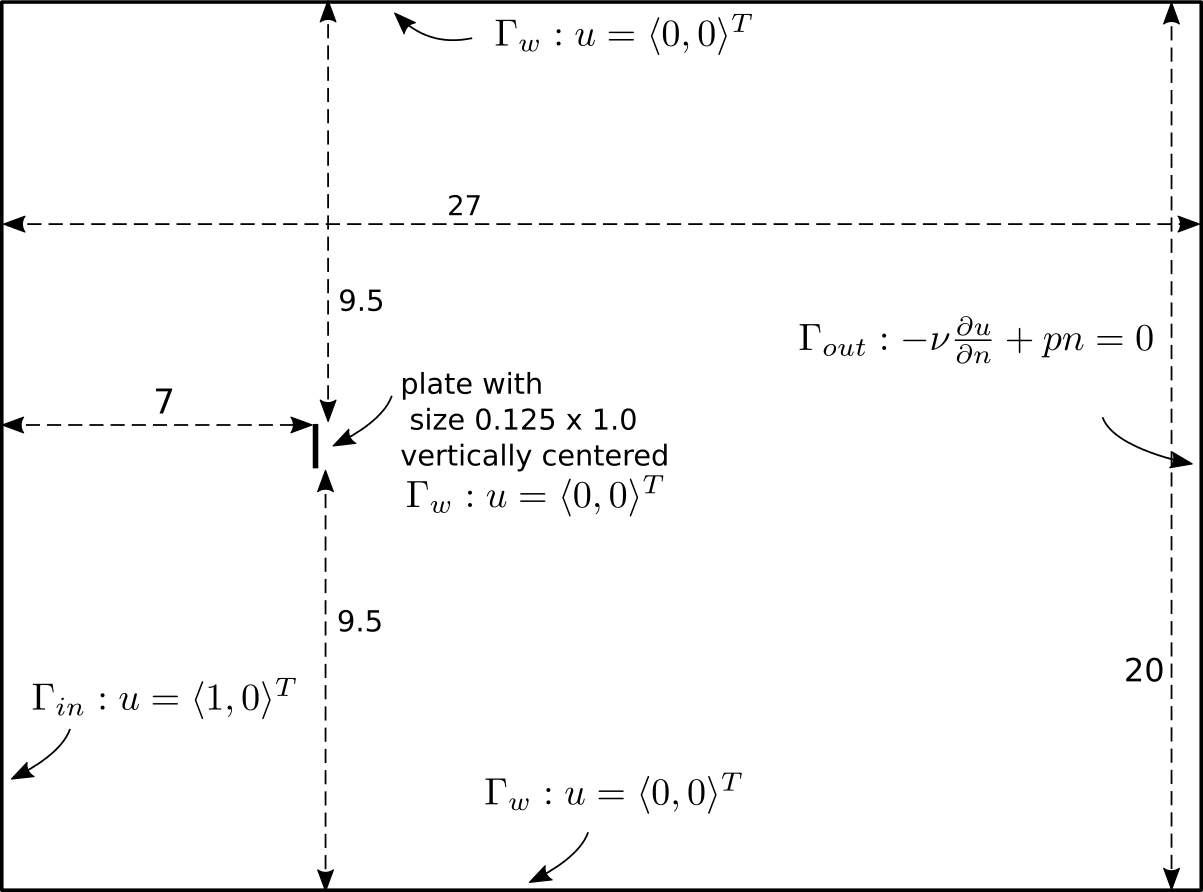}
 \caption{Setup for the flow past a normal flat plate.}
 \label{fig:plate_setup}
 \end{center}
\end{figure}

Our next test is for channel flow past a flat plate with Re=100, following \cite{S07,S13}.  The domain of this test problem is a $[-7,20] \times [-10,10]$ rectangle channel with a $0.125 \times 1$ flat plate placed 7 units into the channel, and vertically centered.  The inflow velocity is set as $\bu_{in}=\langle 1,0 \rangle$, we use a zero-traction outflow, and there is no forcing, ${\bf f} = {\bf 0}$.
No-slip conditions are enforced  on the walls and plate.  A diagram of the test setup is shown in Figure \ref{fig:plate_setup}.

We compute results using the CONV, CONS, SKEW, ROT, and EMAC formulations, with BDF3 time stepping.  The simulations all used BDF3 time stepping, a Delaunay mesh with $(P_2,P_1)$ elements (which provided 58,485 total degrees of freedom) for each simulation.  This is a fairly coarse mesh, and we use it to observe differences between the formulations (since as $h\rightarrow 0$, the formulations will all converge to each other).  The simulations all used the same time step size of $\Delta t=0.02$, and were started from rest, ran until a periodic-in-time state was reached, and then ran for an additional 16 periods.  Periods were determined using the drag coefficient
\begin{align*}
C_d(t^m) & =  \frac{2}{\rho L U_{max}^2} \int_{S} \left( \rho \nu \frac{\partial u_{t_S}(t^m)}{\partial n} n_y - p_h^m n_x \right) \ ds.
\end{align*}
Here, $S$ is the plate, $\bn=\langle n_x,n_y \rangle$ is the outward normal vector, $\bu_{t_S}(t^m)$ is the tangential velocity of $u_h^m$, the density $\rho=1$, the max velocity at the inlet $U_{max}=1$, and $L=1$ is the length of the plate.

The statistics of interest are the average drag coefficient, and the recirculation point of the time averaged velocity; all averages were taken over the last 16 periods.  Results for these statistics are shown in table \ref{platedata}, along with results from a very fine discretization we obtained using the deal.II software \cite{dealII84} and $(Q_2,Q_1)$ elements with the convective formulation and BDF2, using $\Delta t$=0.005 and 4,019,895 total degrees of freedom (for which we {\color{black} assume is a convergent result, since it was very similar to results computed} with $\Delta t=0.01$ and about 2 million total degrees of freedom).  For further comparison, we also give results of Saha from \cite{S07,S13}, who used a MAC scheme with 426x162 cells (16x50 grid points on the plate surface), and a typical time step size of 5E-4.

We note first that the ROT and CONS schemes did not run to completion: the ROT simulation became unstable around T=25, and before T=26 the energy grows to 1E+100; similarly, the CONS scheme gives energy blowup at about T=78.  The EMAC solution's average drag most closely matches that of the very fine discretization, and is significantly closer than that of the CONV and SKEW solutions.  For the recirculation point, the CONV, EMAC, and SKEW formulations give results with similar accuracy.

\begin{table}[h!]
\begin{center}
\begin{tabular}{|l|l|l|l|}
 \hline
Formulation & Re    &  Average $C_d$ & Recirculation point  \\ \hline
%
CONV & 100   &  2.5434 & 1.1577 \\ \hline
EMAC   & 100   &  2.6598 & 1.1648 \\ \hline
SKEW & 100   &  2.5903 & 1.1565 \\ \hline
ROT & 100 & \multicolumn{2}{|l|}{failed: energy blows up at T=25}\\ \hline
CONS & 100 & \multicolumn{2}{|l|}{failed: energy blows up at T=78}\\ \hline\hline
Very fine discretization  & 100   & 2.6454 & 1.1373\\ \hline
Saha \cite{S13} & 100 & 2.43 & 1.11 \\ \hline
Saha \cite{S07} & 100 & 2.60 & (not given) \\ \hline
%
\end{tabular}
\caption{\label{platedata} Shown above are the average drag coefficient and x-coordinate of the recirculation point for simulations of flow past a flat plate with varying formulations, together with reference values from a DNS and from \cite{S07,S13}.}
\end{center}
\end{table}

\subsection{Channel flow past a forward-backward facing step}

Our next experiment concerns flow past a forward-backward facing step.  The domain is a $40\times 10$ rectangle used to represent the channel, and a $1\times 1$ `step' placed at the bottom of the channel, 5 units in.  The boundary conditions are no-slip on all the walls and step, zero-traction at the outflow, and a
constant-in-time parabolic inflow with max inlet velocity of 1.  The initial condition is a parabolic profile across the
channel, and there is no forcing, ${\bf f}={\bf 0}$.  We set the viscosity $\nu=1/600$, and for this setup the correct behavior is for eddies to form behind the step, then detach, move down the channel, and then for new eddies to form, and the cycle repeats \cite{G89,LMNR08}.  The tests are run to an end time of $t=60$.

\begin{figure}[!h]
\begin{center}
{\it CONV} (coarse) \hspace{1.8in} {\it CONV} (finer) \\
\includegraphics[width=.49\textwidth,height=0.2\textwidth, viewport=60 10 500 200, clip]{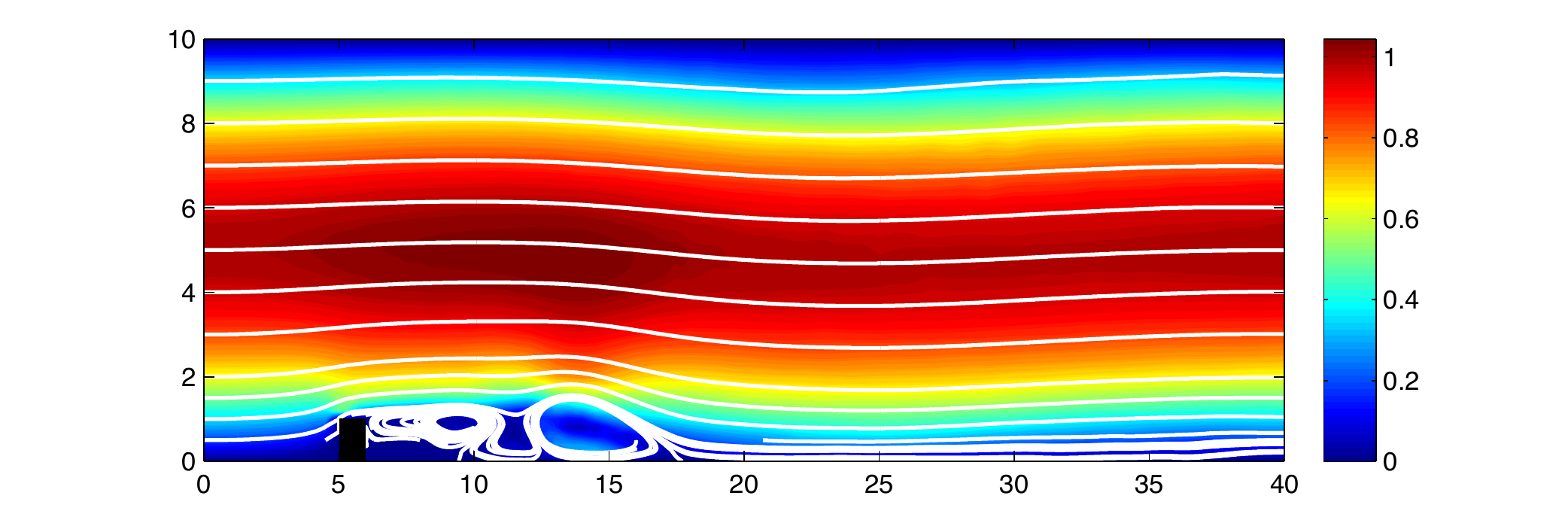}
\includegraphics[width=.49\textwidth,height=0.2\textwidth, viewport=60 10 500 200, clip]{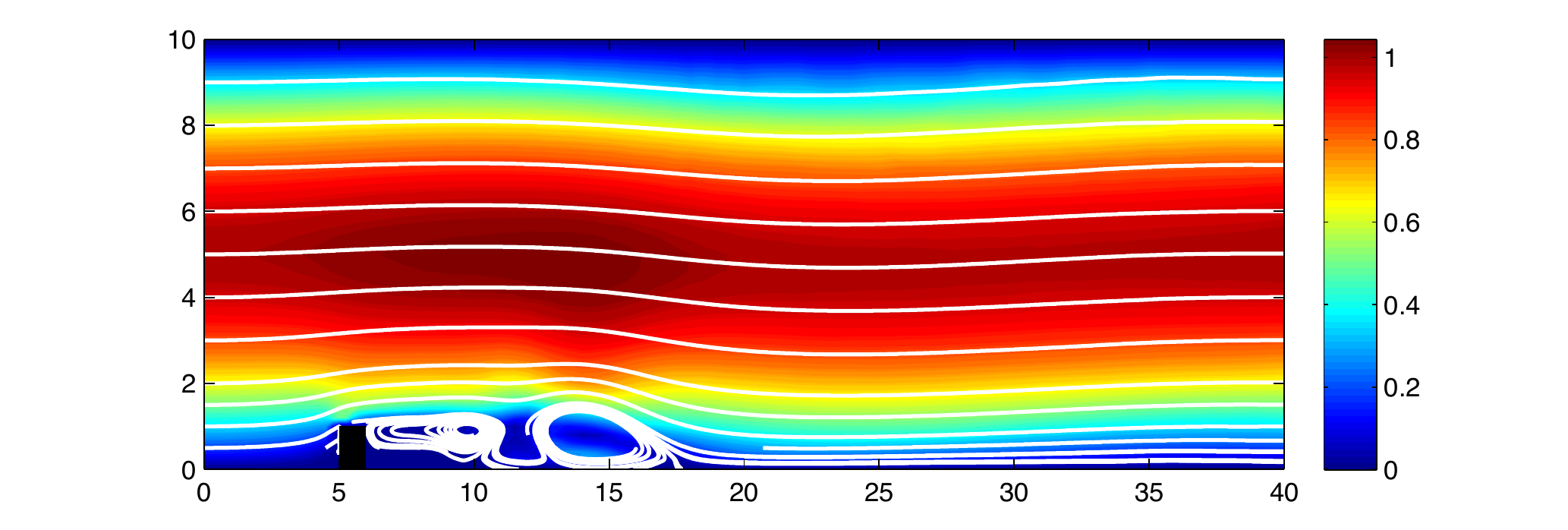} \\
{\it EMAC} (coarse) \hspace{1.8in} {\it EMAC} (finer) \\
\includegraphics[width=.49\textwidth,height=0.2\textwidth, viewport=60 10 500 200, clip]{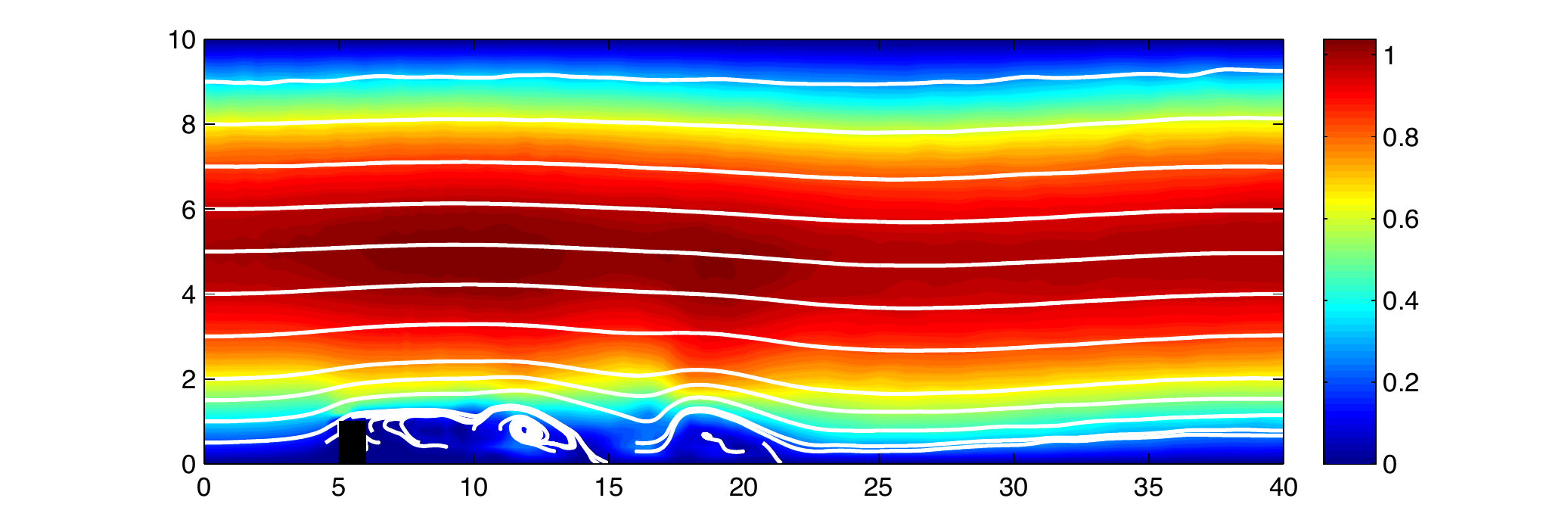}
\includegraphics[width=.49\textwidth,height=0.2\textwidth, viewport=60 10 500 200, clip]{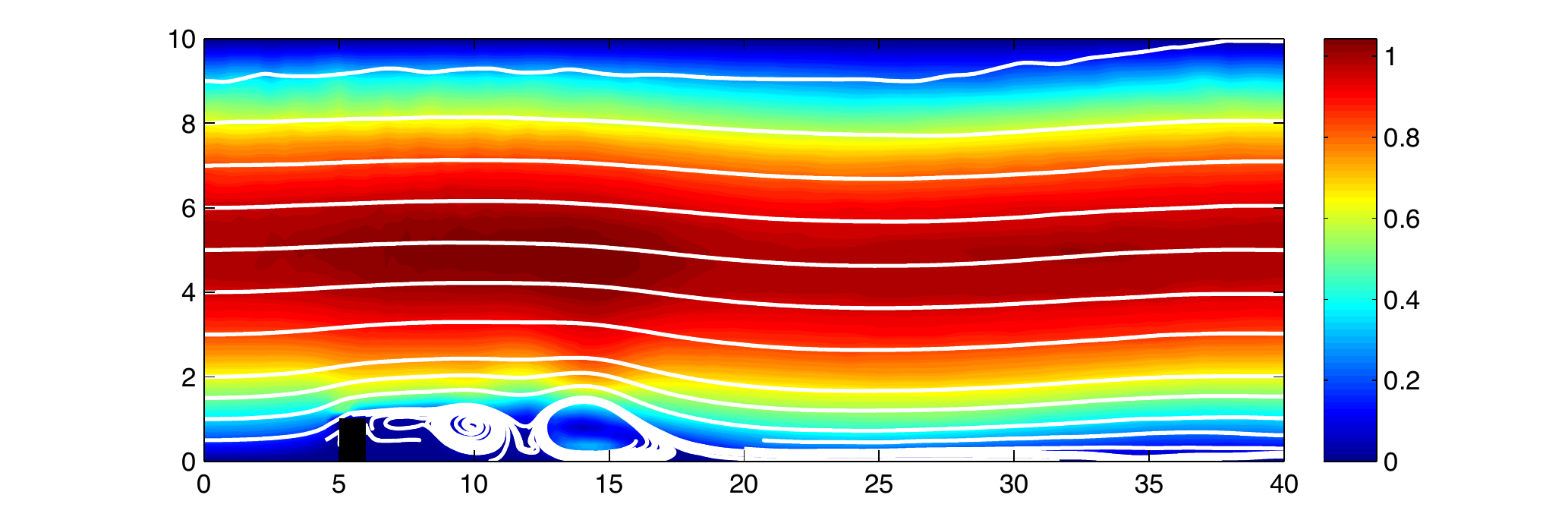} \\
{\it SKEW} (coarse) \hspace{1.8in} {\it SKEW} (finer) \\
\includegraphics[width=.49\textwidth,height=0.2\textwidth, viewport=60 10 500 200, clip]{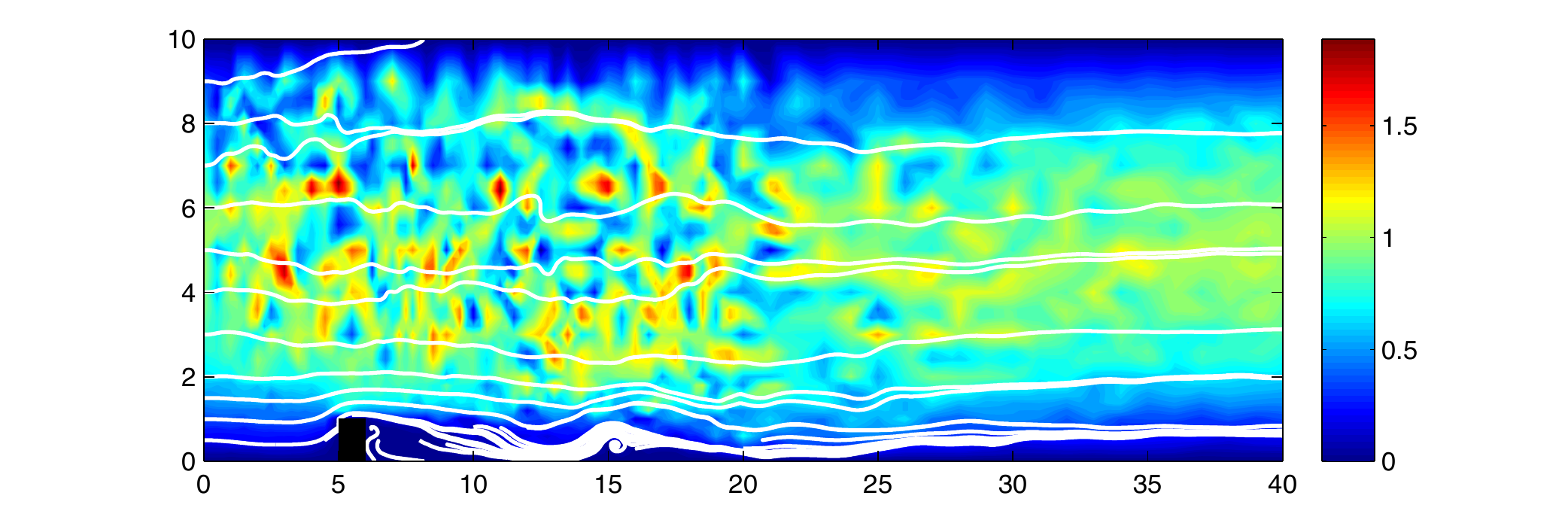}
\includegraphics[width=.49\textwidth,height=0.2\textwidth, viewport=60 10 500 200, clip]{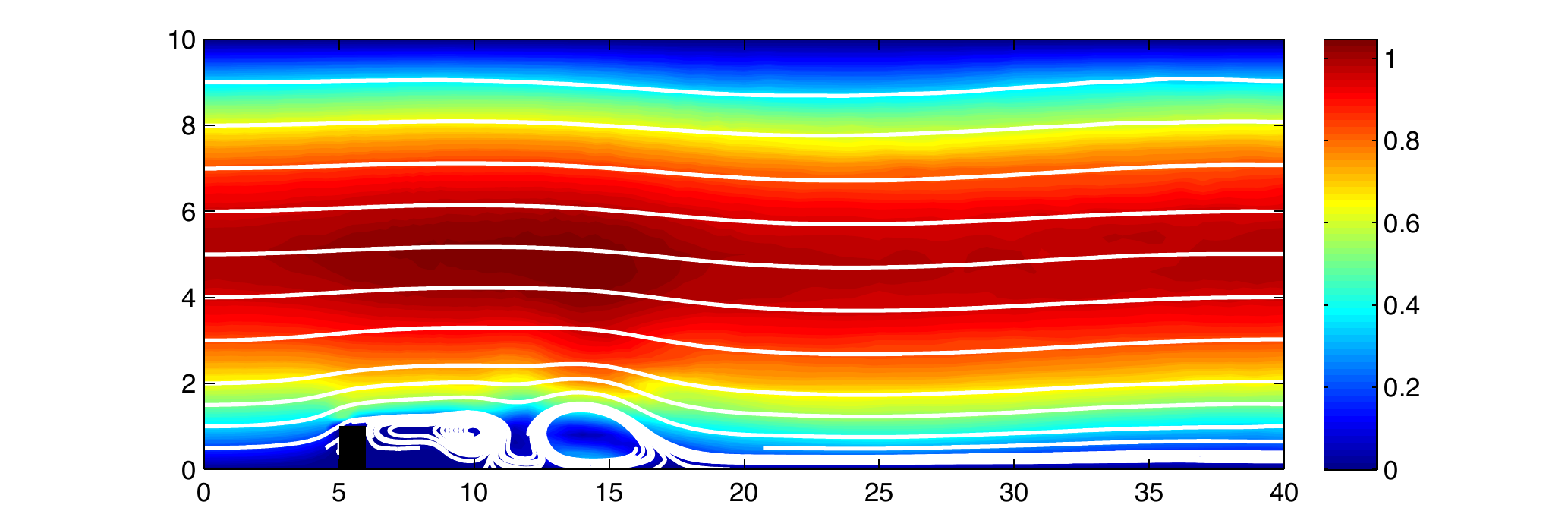} \\
{\it ROT} (coarse) \hspace{1.8in} {\it ROT} (finer) \\
\includegraphics[width=.49\textwidth,height=0.2\textwidth, viewport=60 10 500 200, clip]{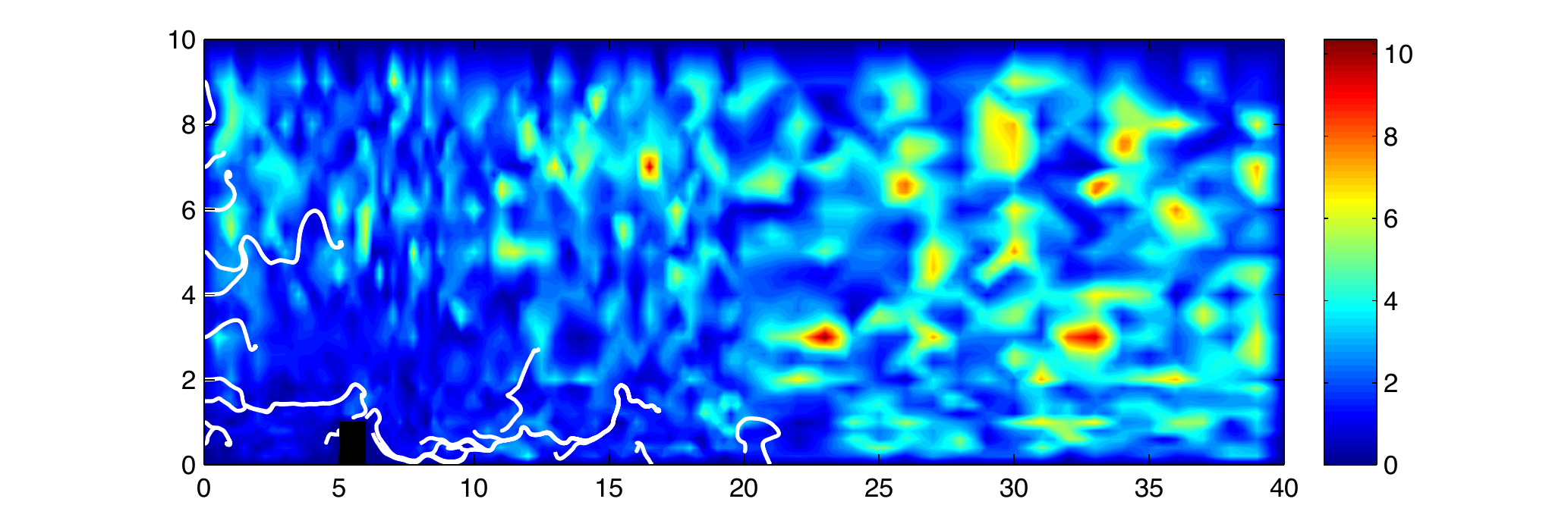}
\includegraphics[width=.49\textwidth,height=0.2\textwidth, viewport=60 10 500 200, clip]{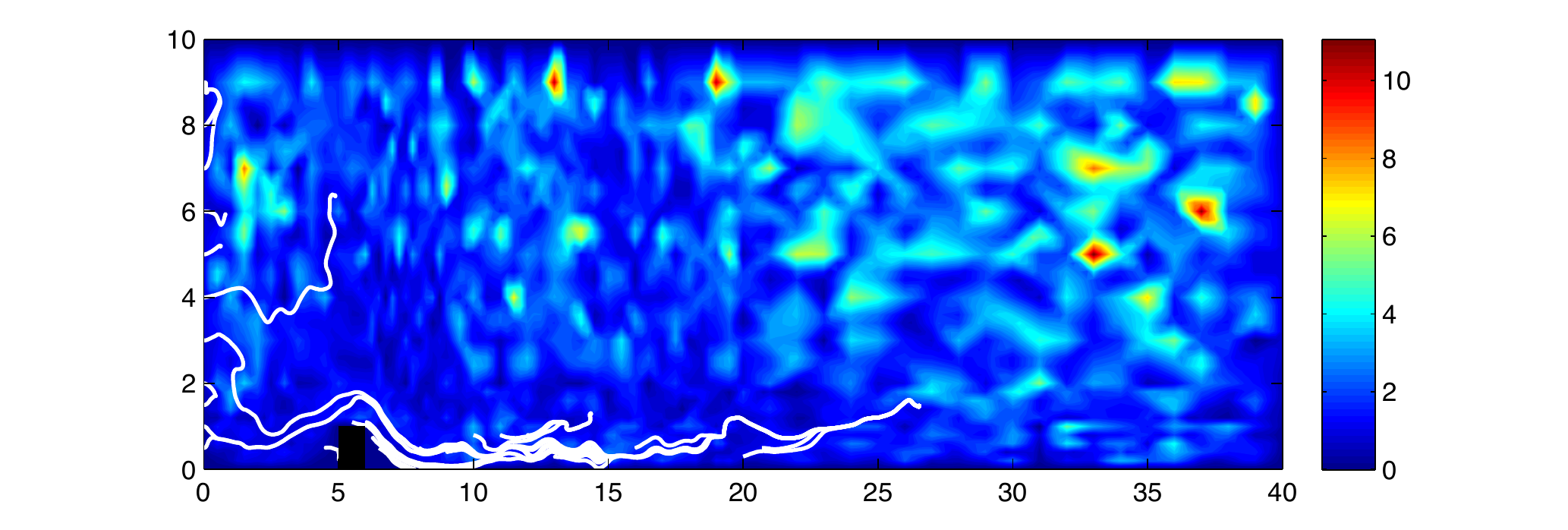} \\
\end{center}
\caption{\label{step2}
Shown above are plots of streamlines over speed contours for the T=60 velocity solutions for flow past a forward-backward step for the different formulations, with the coarser mesh solutions on the left side and the finer mesh solutions on the right.  We note there is no plot for CONS since this simulation failed (energy blowup before T=60).}
\end{figure}

We ran simulations using the five formulations and $(P_2,P_1)$ elements on a coarse mesh that provided 10,908 degrees of freedom, and a finer mesh that provided 19,671 degrees of freedom.  Crank-Nicolson time stepping was used with a time step size of $\Delta t=0.01$.  The CONS scheme failed to finish on both meshes; the energy (rather suddenly) blew up to infinity near t=26 on the coarse mesh, and t=29 on the fine mesh.  The other four schemes all remained stable up to t=60, although with varying accuracies.  Results are shown in Figure \ref{step2} as streamlines over speed contours.  The plots on the right side for EMAC, CONV and SKEW all essentially match the solutions in the literature \cite{LMNR08}, but the ROT solution is quite poor.  On the coarser mesh, the CONV approximation is the best, EMAC is a little worse at resolving the eddies behind the step, and the SKEW and ROT solutions are quite poor.

\subsection{Driven Cavity for Re=10,000}{
\color{black}

In this experiment we consider the well-known 2D flow in a lid-driven square  cavity, see \cite{erturk} for example. We focus on a high Reynolds number of 10,000, which is above the first Hopf bifurcation for Re$\approx$8,000, cf. \cite{auteri2002numerical,bruneau20062dcavity}. Therefore, one expects a periodic flow pattern.

\begin{figure}[!h]
\begin{center}
\includegraphics[width=0.4\textwidth]{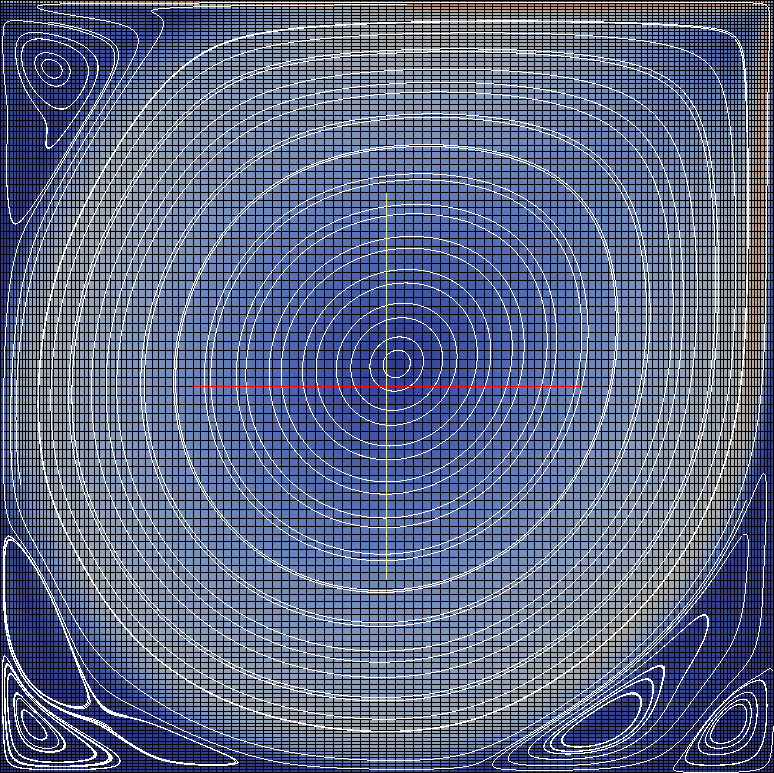}
\includegraphics[width=0.45\textwidth]{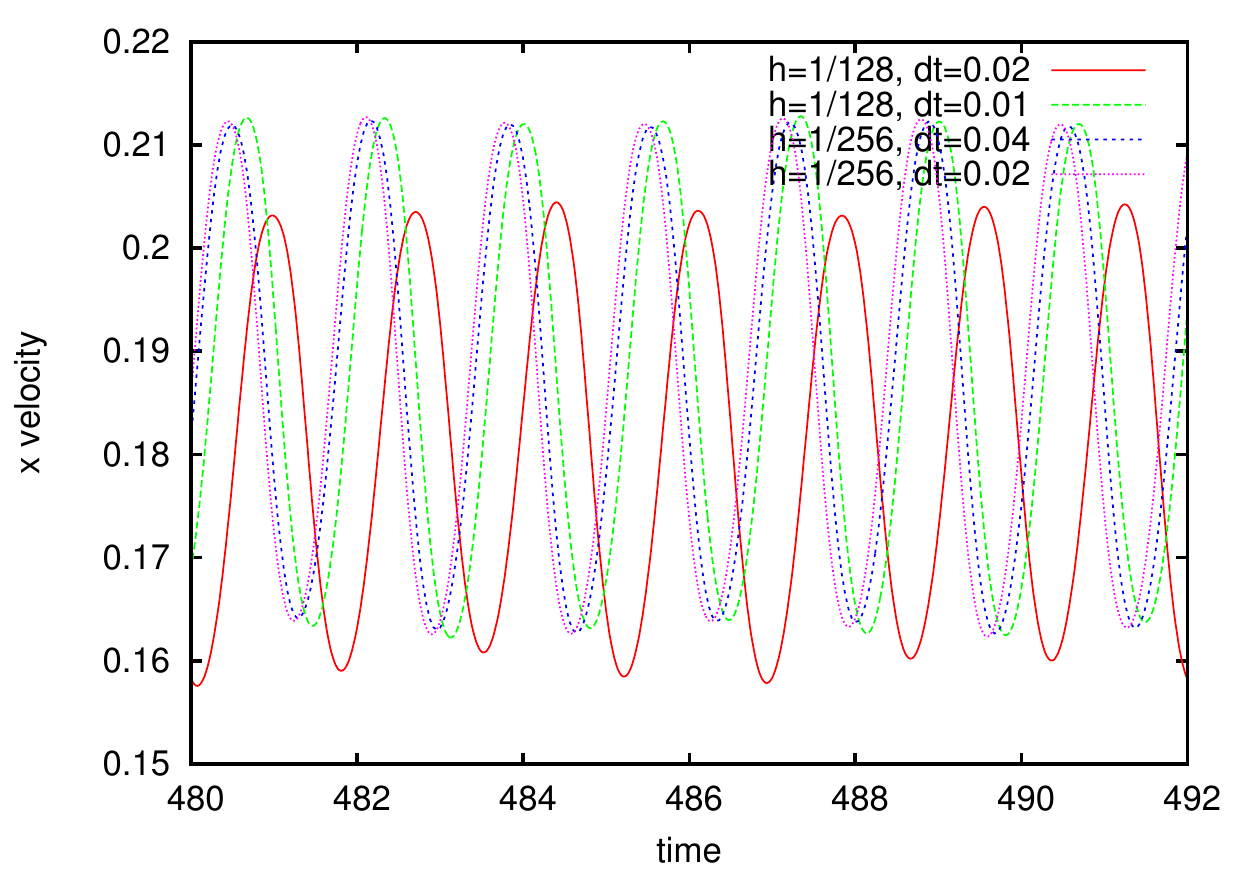}

\end{center}
\caption{2d Cavity with Re=10,000. Left: adapted mesh and velocity streamlines. Right: Plot of x velocity at evaluation point (2/16, 13/16) in periodic regime.}
\label{fig:cavity}
\end{figure}

\begin{figure}[!h]
\begin{center}
\includegraphics[width=0.45\textwidth]{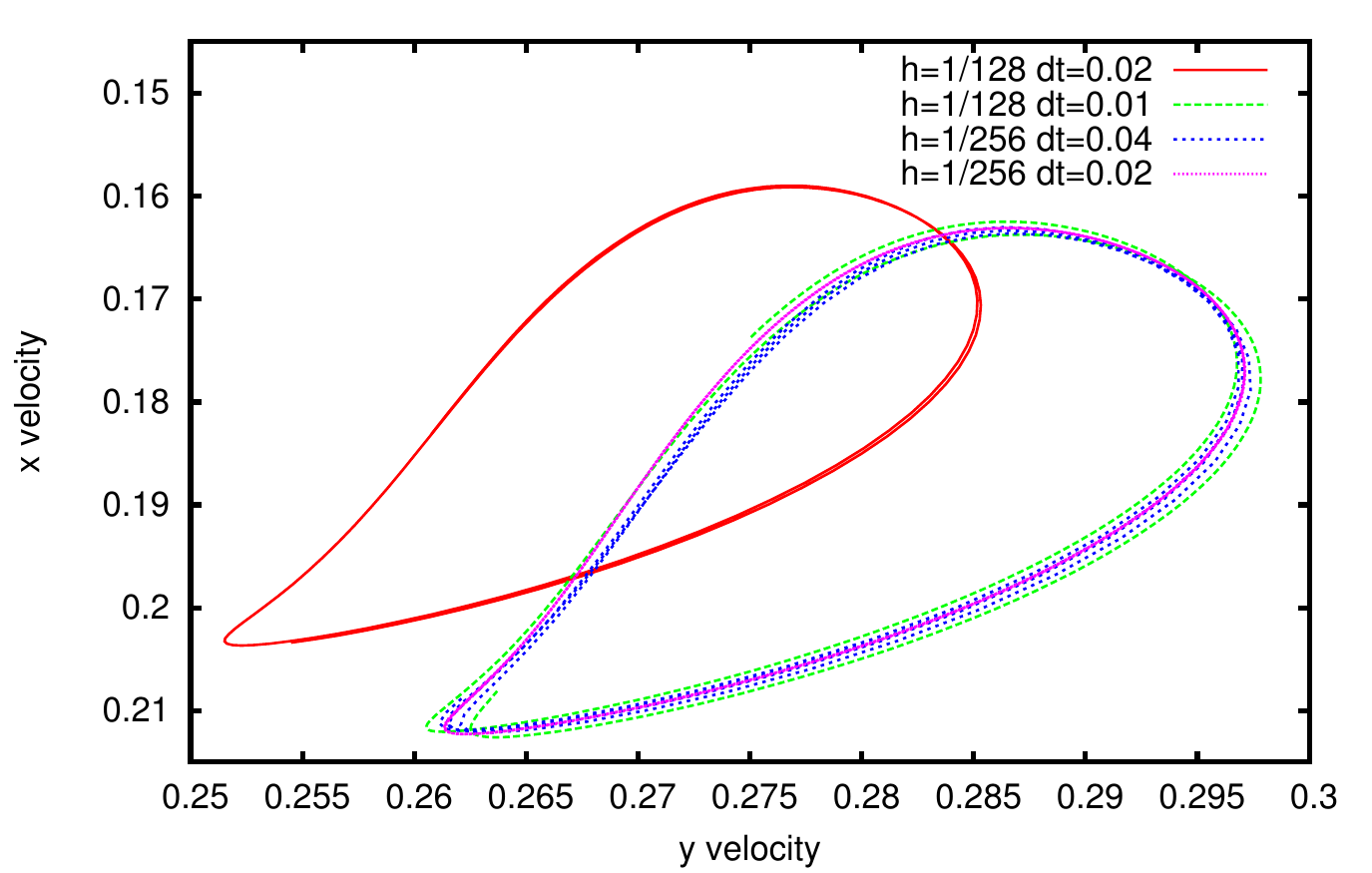}
\includegraphics[width=0.45\textwidth]{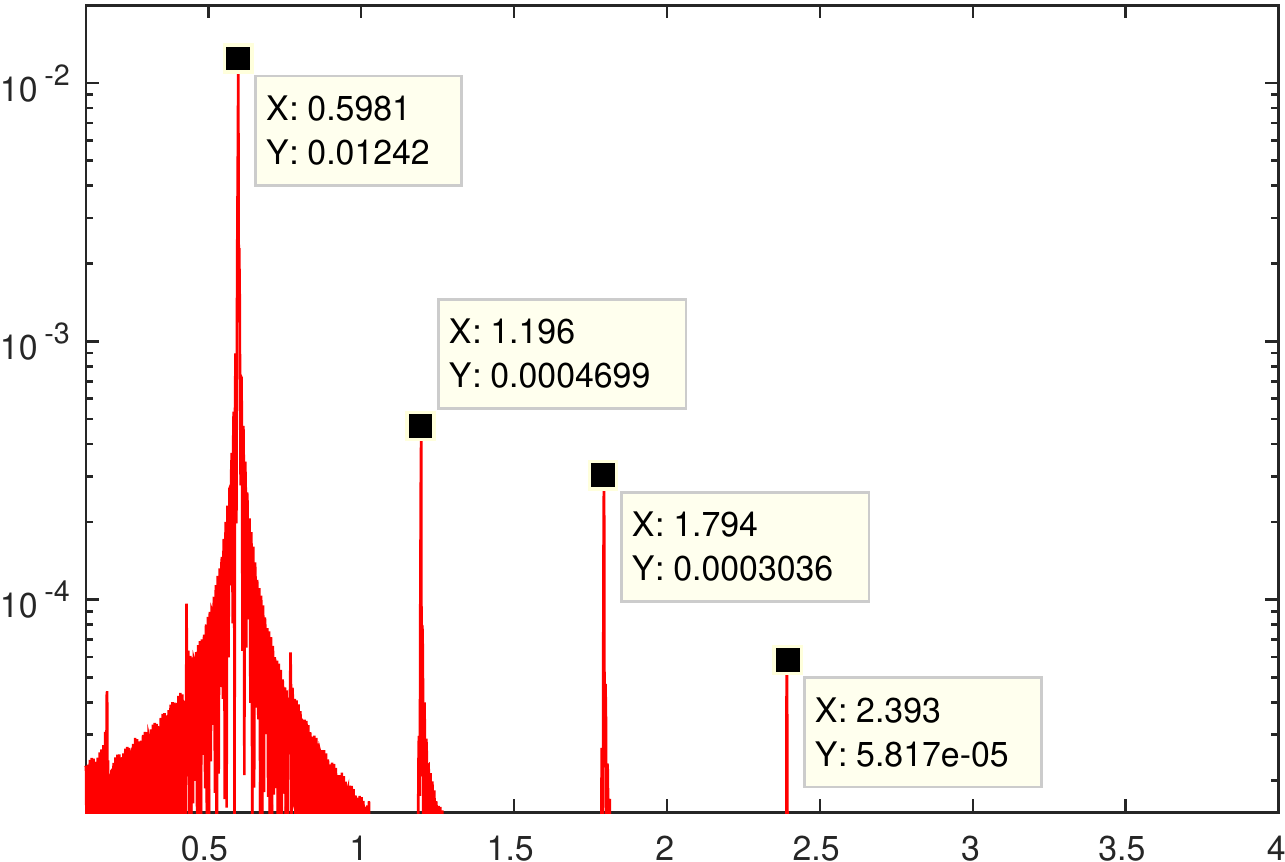}
\end{center}
\caption{2d Cavity with Re=10,000. Left: Phase plot of velocity components at evaluation point (2/16, 13/16) in periodic regime. Right: FFT of velocity component showing the main frequencies.}
\label{fig:cavity2}
\end{figure}

We first compute a stationary solution of the cavity problem for a lower Reynolds (below the Hopf bifurcation). Using this solution as the initial condition,  we proceed by solving the unsteady problem with BDF2 time stepping until we reach a periodic regime. The EMAC scheme is conservative and so only physical dissipation acts to eliminate initial (spurious) frequencies.
Therefore, reaching the periodic state takes a long time.


To resolve the boundary layers and singularities in the corners, we use a wall adapted mesh, as shown in Figure~\ref{fig:cavity}. On this mesh the EMAC formulation was discretized with $(Q_2,Q_1)$ finite elements. As in \cite{bruneau20062dcavity} we evaluate the velocity at the point (2/16, 13/16) {\color{black}near the inflow corner} and plot the x component over time, see Figure~\ref{fig:cavity}, and phase plots and Fourier transforms, see Figure~\ref{fig:cavity2}. Our statistics demonstrate mesh convergence and  show that we resolve the problem using a timestep of 0.02 and a mesh size of 1/256 (592,387 DoFs). The computed  statistics are in good agreement with \cite{bruneau20062dcavity} even with a much coarser mesh than the one  used in that paper.
}

\subsection{3D Flow around a Square Cylinder}

As a final test, we include some 3d flow computations of the new EMAC formulation in order to validate
that the scheme generates the correct periodic flow behavior for a flow in a channel past a square cylinder {\color{black}
at Re=100}. The setup of the benchmark problem is taken from \cite{schafer1996benchmark}
(unsteady test case 3D-2Q) and we will compare with results from \cite{schafer1996benchmark} and \cite{Olshanskii2013231}.
{\color{black}The Reynolds number of 100 is close to the critical one, where the transition from equilibrium to unsteady periodic solution takes place
for this problem.   This makes the test challenging for a NSE numerical solver, since the discretization  method should ensure the right balance between inertia
and viscous  diffusion to produce solutions which are numerically stable and, at the same time, the development of unsteady behavior is not suppressed~\cite{Olshanskii2013231}.}

\begin{figure}[!h]
\begin{center}

\includegraphics[width=0.85\textwidth]{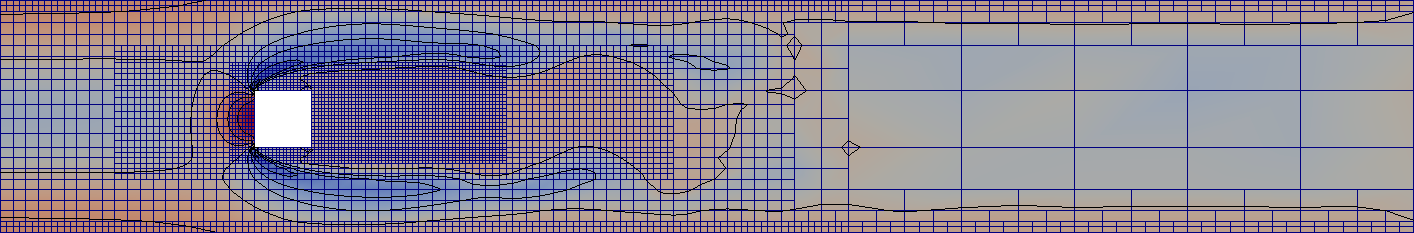}

\vspace{1em}

\includegraphics[width=0.85\textwidth]{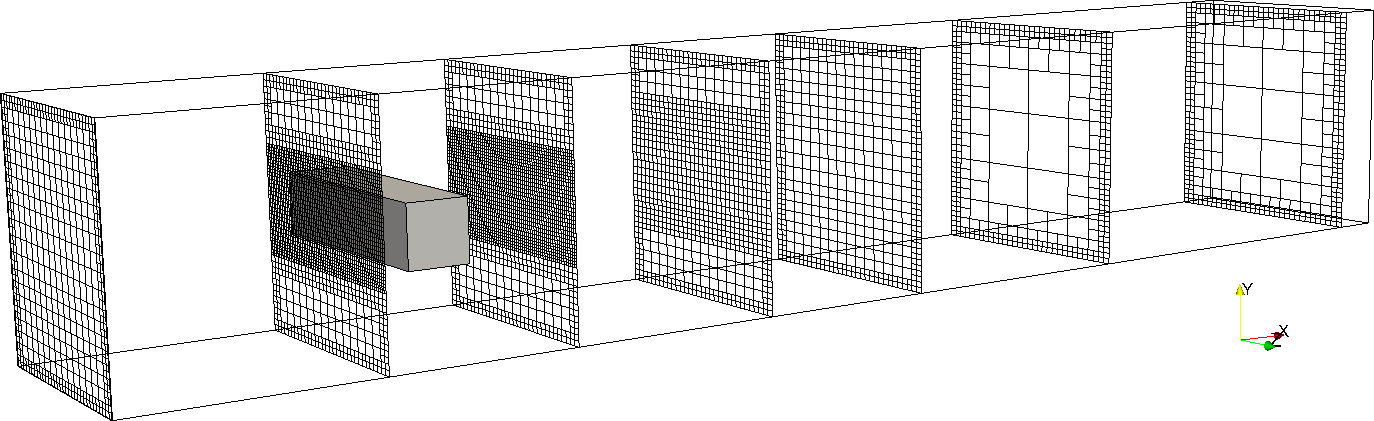}

\includegraphics[width=0.45\textwidth]{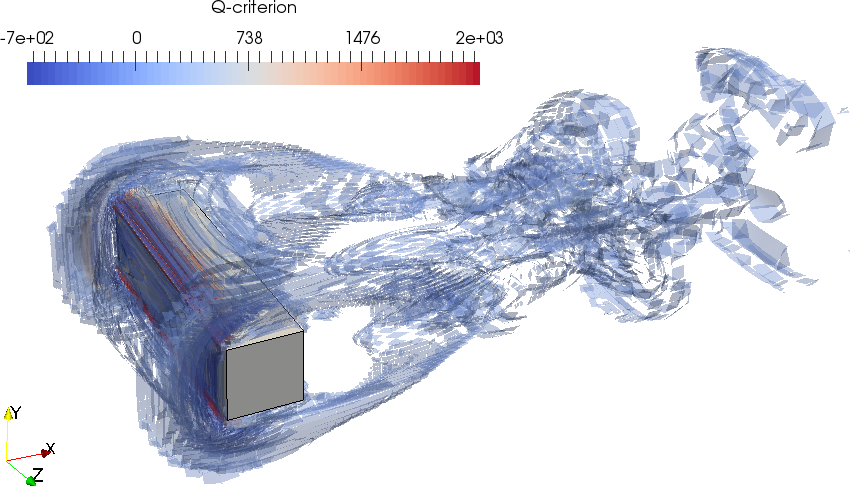}
\includegraphics[width=0.4\textwidth]{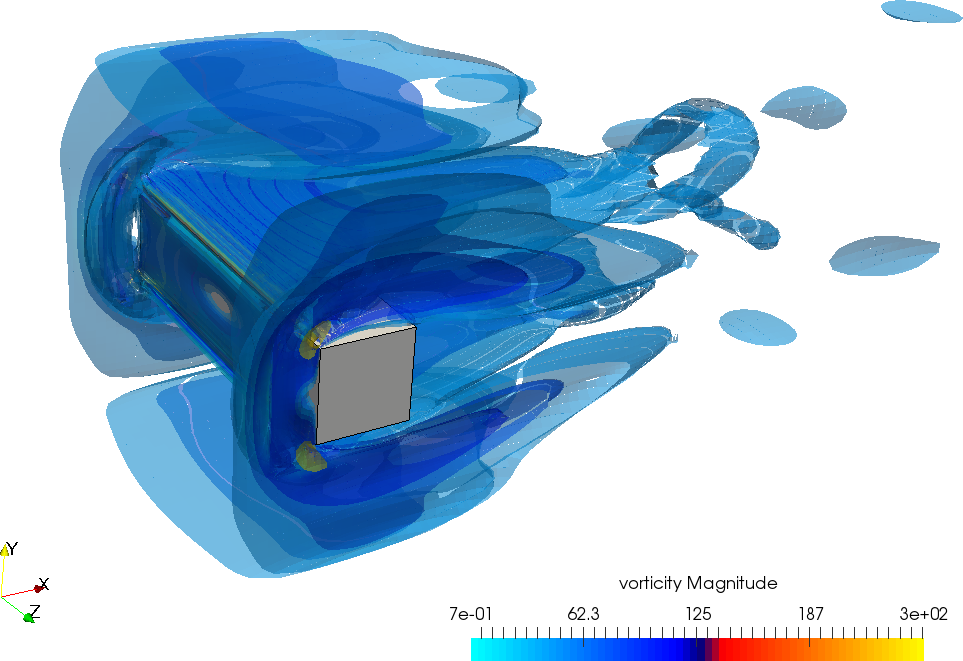}

\end{center}
\caption{3D flow around a square cylinder. Top: mesh through the midplane, coloring by pressure. Middle: slices of the
mesh on various downstream planes. Bottom: q-criterion and vorticity magnitude contours at t=12.0 that show the 3d structure of the solution.}
\label{fig:flow3d-mesh}
\end{figure}

The computation is done for the time interval T=[0,13] with step size 0.01 using BDF2
time discretization.
As in the cavity problem, the nonlinear problem of the EMAC scheme in each time step is discretized using  $(Q_2,Q_1)$ finite elements on a fixed mesh, which is displayed
in Figure~\ref{fig:flow3d-mesh}. The mesh has been refined manually
to about 240k cells resulting in 6.4 million DoFs.
The linear systems
are again solved like in the example before using a grad-div parameter of $\gamma=0.1$.

\begin{table}
 \centering
 \begin{tabular}{|l|l|l|l|l|} \hline
   & max drag & max lift & Strouhal & DoFs \\ \hline
  EMAC results & 4.890 & 0.0271 & 0.351 & 6.4 mill \\
  \cite{Olshanskii2013231} & 4.484 & 0.0316 & 0.307 & 17 mill \\
  \cite{schafer1996benchmark} & 4.32-4.67 & 0.015-0.05 & 0.27-0.35 & Up to 6 mill \\ \hline
  \end{tabular}
\caption{Table of reference values of the flow around the 3d square cylinder.}
\label{tab:3dflow-results}
\end{table}

\begin{figure}[!h]
\begin{center}
\includegraphics[width=0.45\textwidth]{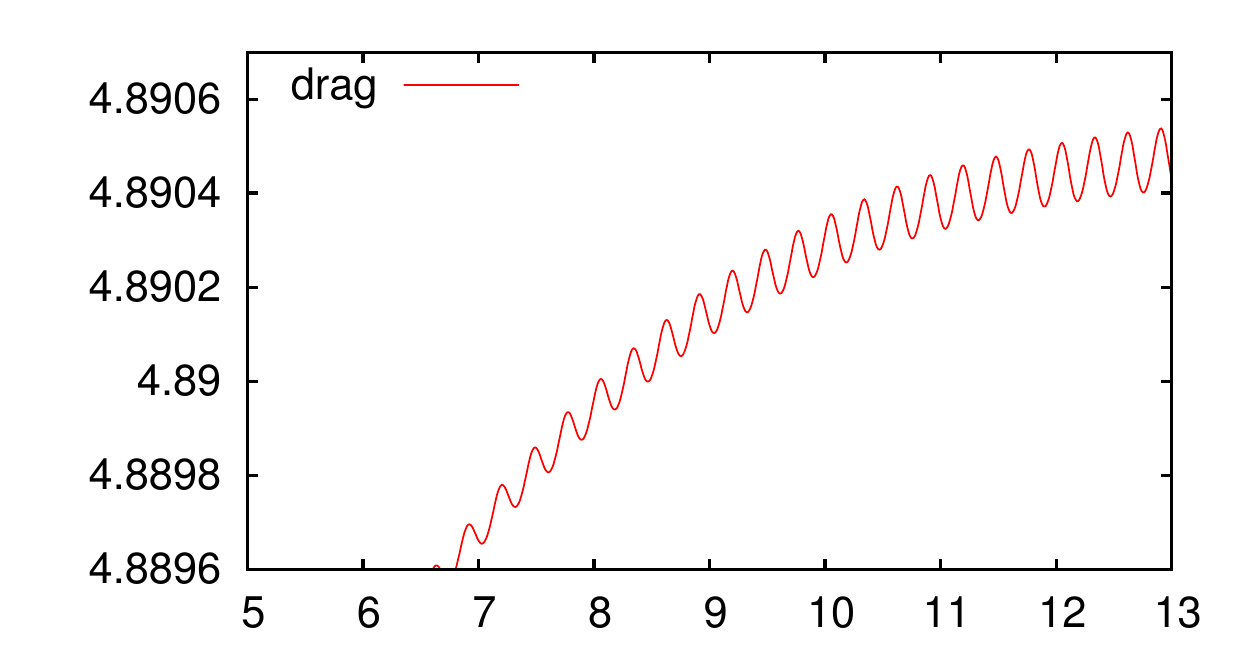}
\includegraphics[width=0.45\textwidth]{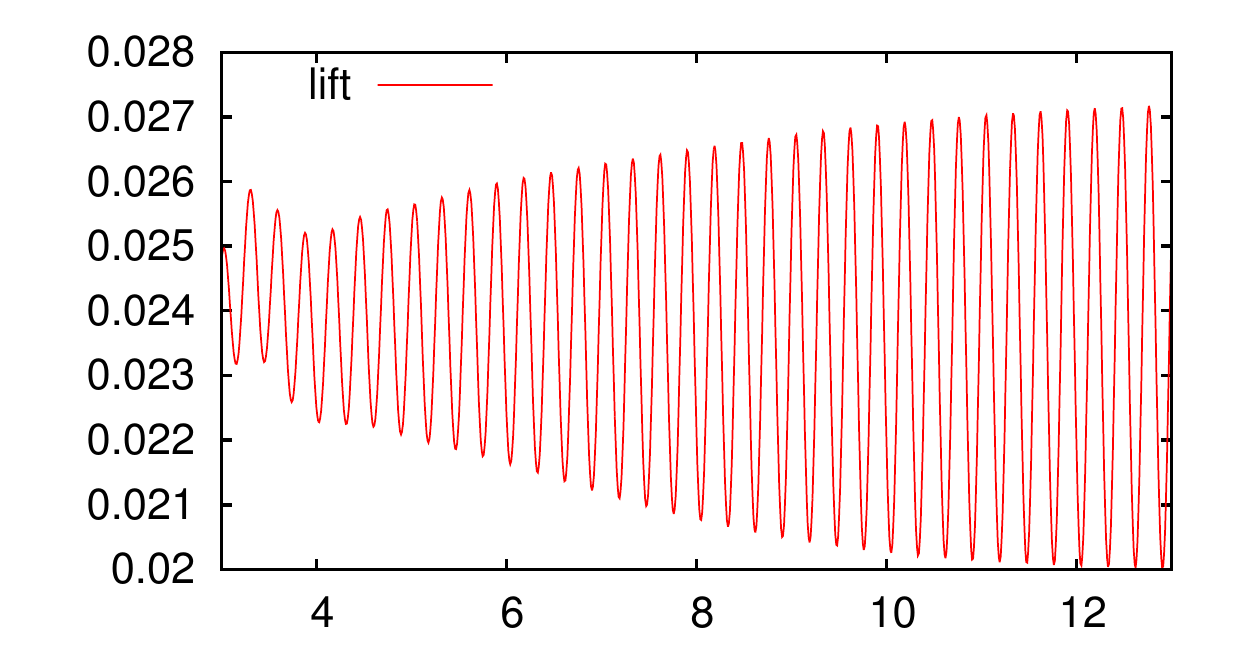}
\end{center}
\caption{Lift and drag coefficients of the flow around the 3d square cylinder plotted over time T.}
\label{fig:flow3d-liftdrag}
\end{figure}

The results in Table~\ref{tab:3dflow-results} show good agreement with the reference values from the
literature even though the mesh is relatively coarse.  {\color{black}Note that \cite{schafer1996benchmark} does not give reference intervals for this problem, and we simply show the maximum  and minimum values of   lift, drag, and the Strouhal numbers for several DNS results  included in \cite{schafer1996benchmark}. However, these intervals can be not very accurate.}
The visualizations in Figure~\ref{fig:flow3d-mesh} show that 3d flow structures develop behind the cylinder as expected.

\section{Conclusions and Future Directions}
We have developed a new discrete formulation for incompressible Navier-Stokes equations, named the {\it EMA-conserving} (EMAC) formulation herein,
which conserves energy, momentum, angular momentum, and appropriately defined vorticity, helicity, and enstrophy, when the solenoidal constraint on the velocity is enforced only weakly.  {\color{black}Moreover, we show that none of the commonly used {\it convective}, {\it conservative}, {\it rotational}, and {\it skew-symmetric} formulations conserve each of energy, momentum, and angular momentum.  
These properties together with the form of non-linear and potential terms are summarized in the table below.
\begin{table}[h!]
\small
\begin{tabular}{rl|c|c|c|c}
name& $NL(\bu)$  & potential term                                  &  Energy & Momentum & Ang. Moment. \\ \hline
convective: & $\bu \cdot\nabla \bu$ &$p$ (kinematic)             &&&\\
skew-symm.: & $\bu \cdot\nabla \bu +  \frac12 (\div \bu)\bu$&$p$ (kinematic) &\textbf{+}&&\\
rotational: & $(\nabla\times\bu)\times \bu$&$p+\frac12|\bu|^2$ (Bernoulli)       &\textbf{+}&&\\
conservative: & $\nabla\cdot (\bu{\color{black}\otimes}\bu)$&$p$ (kinematic)&&\textbf{+}&\textbf{+}\\
EMAC: &$2\bD(\bu)\bu + (\div \bu)\bu$&$p-\frac12|\bu|^2$ (no name)  &\textbf{+}&\textbf{+}&\textbf{+}
\end{tabular}
\end{table}
Results of several numerical experiments have been provided which verify the discrete conservation properties of the EMAC scheme, and also show that it performs at least as good, or better, than the commonly used formulations.}

Aside from further testing, one important future direction is to consider more efficient treatments of the {\it EMA-conserving} formulation.  That is, in this initial study, we consider schemes that solve the nonlinear problem at each timestep.  However, it is typical with the more commonly used formulations to linearize the nonlinear term at each time step by approximating one of the velocities using previous time step solutions; such schemes need only one linear solve per time step, whereas schemes that resolve the full nonlinear problem with Newton's method often require two or three.


\section*{Acknowledgements}
{\color{black}
Clemson University is acknowledged for generous allotment of compute time on Palmetto cluster.
}

\bibliographystyle{plain}
\bibliography{graddiv}

\end{document}